\documentclass{article}

\usepackage{amsfonts}
\usepackage{float}
\usepackage{amsmath}
\usepackage{amssymb}
\usepackage{amsthm}
\usepackage{graphicx}
\usepackage{multirow}
\usepackage{color}
\usepackage{ulem}
\usepackage{setspace}
\usepackage{caption}
\usepackage{subcaption}
\usepackage{rotating}
\usepackage{threeparttable}
\usepackage{url}
\usepackage{framed}

\newcommand{\mbf}[1]{\mbox{\boldmath $#1$}}

\def\bA{\mathbf{A}}

\def\bR{\mathbf{R}}

\def\bW{\mathbf{W}}
\def\bX{\mathbf{X}}

\def\bD{{\mathbf D}}

{\catcode `\@=11 \global\let\AddToReset=\@addtoreset}
\AddToReset{equation}{section}

\AddToReset{Theorem}{section}

\newtheorem{lem}{Lemma}[section]
\newtheorem{thm}{Theorem}[section]

\theoremstyle{remark}
\newtheorem{rem}{Remark}[section]

\newcommand{\cB}{{\cal B}}

\newcommand{\cD}{{\cal D}}

\newcommand{\cL}{{\cal L}}

\newcommand{\cN}{{\cal N}}

\newcommand{\cQ}{{\cal Q}}

\newcommand{\cW}{{\cal W}}

\def\bc{\begin{center}}
\def\bd{\begin{description}}
\def\be{\begin{enumerate}}
\def\ec{\end{center}}
\def\ed{\end{description}}

\def\ee{\end{enumerate}}
\def\ben{\begin{equation}}
\def\benn{\begin{equation*}}
\def\een{\end{equation}}
\def\eenn{\end{equation*}}
\def\benr{\begin{eqnarray}}
\def\eenr{\end{eqnarray}}
\def\benrr{\begin{eqnarray*}}
\def\eenrr{\end{eqnarray*}}

\def\b{\beta}

\def\edt{\end{document}}

\def\G{\boldsymbol{\Gamma}}

\def\iny{\infty}

\def\lel{\label}

\def\noi{\noindent}

\def\ra{\rightarrow}

\def\si{\sigma}

\def\stt{\sum_{t=1}^n}
\def\sti{\sum_{i=1}^n}

\def\stj{\sum_{j=1}^J}
\def\stk{\sum_{k=1}^n}
\def\stl{\sum_{l=1}^n}
\def\sth{\sum_{h=1}^n}

\def\vep{\varepsilon}

\def\wh{\widehat}

\def\wt{\widetilde}

\def\R{{\mathbb R}}

\def\bx{\mbf{x}}

\def\B{\boldsymbol{\beta}}
\def\W{\boldsymbol{\cW}}
\def\S{\boldsymbol{\Sigma}}

\def\L{\boldsymbol{\Lambda}}

\def\tbf{\textbf}
\def\tecr{\textcolor{red}}

\newcounter{Qcounter}
\newcounter{show}
\title{Robust and efficient estimation of binary regression parameters via solving a $L_{2}$ optimization problem}
\author{Jiwoong Kim\\
Department of Mathematics \& Statistics\\
University of South Florida}

\begin{document}
\maketitle

\begin{abstract}
The statistical literature is replete with research on binary regression, and it is not difficult to find statistical software or packages that provide statistical inference for it. This paper proposes a novel method for parameter estimation in the binary regression model. Upon obtaining estimators by the proposed method, this study rigorously investigates their asymptotic properties and computational aspects. To validate the proposed method, this study compares it with other well-celebrated methods. Empirical studies demonstrate that the proposed method compares favorably with other methods.
\end{abstract}

\section{Introduction}
After Cram\'{e}r (1928) and von Mises (1928) proposed the Cramer-von Mises (CvM) criterion, it has been extensively used for the goodness-of-fit (GOF) test: see, e.g, Anderson and Darling (1952) and Anderson (1962). Since then, it has been referred to as the CvM test, rather than the CvM criterion. The application of the CvM criterion was not limited to the GOF test; it has been widely adopted by statisticians for both parametric and nonparametric estimation methods. One of the well-known estimation methods, in which the CvM criterion has left its seminal mark in the statistical literature, is minimum distance (MD) estimation. In the literature on the MD estimation, the CvM-type distance refers to any measure that integrates the difference between two entities that appear parallel to each other. For example, the probability measure in the original CvM criterion has been replaced by many other integrating measures (including Lebesgue and degenerate) for application purposes.

To briefly recap the MD estimation, it has been widely used, with applications ranging from the estimation of a location parameter in a one-sample model to the estimation of parameters in regression and autoregressive models: see, e.g., Koul (2002) and references therein. However, all the applications mentioned above are limited to cases where the sample of random variables follows a continuous distribution: the MD estimation using the CvM-type distance has not been applied to either a sample of discrete random variables or to regression models with discrete response variables. Giving a second thought about why it was not used at all for the discrete cases, one can easily find an answer. Consider, for example, a linear regression model with a discrete response term. The CvM-type $L_{2}$ distance -- the most popular and common distance function in the MD estimation -- integrates the  square of difference between the empirical and modeled distribution functions of the response variable. The existence of the integrating measure embedded in the distance function -- e.g., the Lebesgue measure -- explains why the CvM type distance estimation is intrinsically the least amenable of the MD estimations to being applied to the discrete random variables: it requires the modeled distribution function to be continuous, which can be utterly impossible for discrete random variables.

Considering its asymptotically desirable properties, such as efficiency and robustness, it is indeed regrettable that the MD estimation with the CvM-type distance can not be applied to the parameter estimation of discrete distributions. Hence, its applicability to discrete distributions merits further investigation. Motivated by its desirable properties, Kim (2025) applied the MD methodology to estimate the success probability in the binomial distribution, replacing the continuous integration measure with a discrete measure. Kim (2025) demonstrated the asymptotic normality of the resulting MD estimator; furthermore, he empirically showed that the MD estimator is more robust than the maximum likelihood (ML) estimator in the presence of outliers. He also compared it with other well-celebrated robust estimators, including the E-estimator and the entropy-based divergence estimator: for more details on these estimators, see Jones et al. (2001), Ruckstuhl and Welsh (2001), Fujisawa and Eguchi (2008), and Kawashima and Fujisawa (2017). Motivated by the findings in Kim (2025), we extend the domain of application of MD estimation from one sample of discrete random variables to a regression setup, viz., a binary regression model.

The rest of this article is organized as follows. Section \ref{sec:MDE} introduces a modified version of the MD estimation method tailored to the current binary regression setup, proposes the MD estimator of the regression parameters, and presents the asymptotic properties of the proposed estimators, which are the main results of this study. Section \ref{Sec:emp_studies} compares the proposed estimation method with other methods for binary regression, including the generalized linear model (GLM) and its variants, and presents the results of the comparison based on empirical findings that comprise simulation studies and real data examples. The statistical package \texttt{jwglm}, used in the empirical studies presented in this article, is available on \url{https://github.com/jwboys26}. Section \ref{Sec:conclusion} provides a brief recap of this study and concludes the article. For a real-valued function $f$ mapping $\mathbb{R}$ to $\mathbb{R}$, $\|f\|_{\iny}$ will denote its supremum norm, that is, $\|f\|_{\iny}:=\sup_{x\in\mathbb{R}}|f(x)|$. To denote non-scalar values, boldface will be used. For example, an italicized, small boldface denotes a vector, while an upright, capital boldface denotes a matrix. For a real vector $\mbf{a}$, $\|\mbf{a}\|$ will denote its Euclidean norm. For a function $f$, $f'$ will denote its first derivative, while the prime will denote a transpose when it is used with a vector or a matrix: e.g., $\mathbf{M}'$ will denote the transpose of a matrix \tbf{M}, whereas $\mbf{x}'$ will denote the transpose of a vector $\mbf{x}$.

\section{Minimum distance estimation}\lel{sec:MDE}
\subsection{More literature review of the MD estimation}\lel{sec:Literature}
After Wolfowitz (1953, 1957) published a series of research articles on the MD estimation, it attracted the attention of many researchers. Among them, Parr and Schucany (1980), Millar (1981, 1982, 1984), and Donoho and Liu (1988a, b) conducted exhaustive research on the MD estimation, wherein they explored the MD estimation with various distances, including the CvM, Kolmogorov-Smirnov, Kuiper, Hellinger, L\'{e}vy, and Prohorov distances. After employing several robustness criteria proposed by Hampel (1968, 1974, 1986), Parr and Schucany (1980) empirically demonstrated that the robustness of the MD estimation can be achieved by using the CvM-type distance. However, the research mentioned above applied the MD methodology only to the one-sample model. Extending the application of the MD estimation to the regression and autoregressive setups, Koul (1970, 1985a,b, 1986, 2002) rigorously investigated asymptotic properties of the MD estimators – obtained from the CvM-type distance – of the parameters of regression and autoregressive models where the error and innovation in the models are independent. Koul (2002) demonstrated that the proper choice of weights for the empirical process embedded in the distance function gives rise to efficient estimators of parameters in the linear model for various error distributions, including normal, logistic, and double exponential distributions: see Chapter 5.6 of Koul (2002) for more details. Afterward, other researchers have conducted further studies on extending the MD estimation using the CvM-type distance. Kim (2020) studied the MD estimation of parameters of a regression model with dependent errors. For the computational aspect, Dhar (1991, 1992) demonstrated that the MD estimators of regression and autoregressive parameters exist under certain conditions. Kim (2018) proposed a novel algorithm that yields the MD estimator much more efficiently than other well-celebrated numerical methods.

\subsection{The MD estimation for binary regression}\lel{sec:ULAQ}
Before commencing the MD estimation method, we will provide a summary of binary regression in this section. The binary regression is a statistical model that constructs a relationship between a binary response variable $y$ and either categorical or continuous predictors $\mbf{x}:=(x_{1},...,x_{J})\in \mathbb{R}^{J}$, where the probability that $y$ is happening can be modeled as
$\mathbb{P}(y=1) = p(\mbf{x}'\B)$, where $p:\mathbb{R}\ra [0,1]$. Note that
$p$ is indeed a distribution function (df), and popular examples of $p$ include; a logistic $p(s):=1/(1+e^{-s})$; a normal $p(s):=\Phi(s)$ where $\Phi$ is a df of a normal random variable; and a Cauchy $p(s)=\pi^{-1}\arctan(s)+0.5$. Among all, the logistic df has been the most popular for $p$ with researchers. Since binary regression using the logistic $p$ has wide currency in the statistical literature, the umbrella term ``logistic regression" is commonly used, rather than binary logistic regression, demonstrating the popularity of the logistic $p$. In this article, the investigation will not be confined to the binary logistic regression unless specified otherwise.

Note that the binary regression model builds a relationship between the probability of a response variable and predictors, which is the fundamental difference from other linear regression models that build a direct relationship in that the response variable itself is employed in other models. Having a root cause in this fundamental difference, the binary regression is more ``computationally-expensive" than other linear models since estimators of the binary regression parameters do not have a closed-form solution and can be obtained only via numerical optimization, such as the gradient descent (GD) method. Consider $n$ pairs of observations $(y_{1}, \mbf{x}_{1}),...,(y_{n}, \mbf{x}_{n})$. Based on the observations, the likelihood for the binary regression is defined as
\benrr
L(\B)&=& \prod_{k=1}^{n}p(\mbf{x}_{k}'\B)^{y_{k}}(1-p(\mbf{x}_{k}'\B))^{1-y_{k}},
\eenrr
or, equivalently, the log-likelihood can be defined as
\ben\lel{eq:log_likelihood}
l(\B) = \stk\left[y_{k}\log p(\mbf{x}_{k}'\B) +(1-y_{k})\log (1-p(\mbf{x}_{k}'\B) ) \right],
\een
where $p$ is a probability df. Then, the maximum likelihood (ML) estimator -- the most popular estimator for the binary regression in the statistical literature -- is defined as
\benn
\hat{\B}_{ML} = \arg\max_{\B\in\R^{J}} l(\B).
\eenn
It is known that the ML estimator exhibits significant bias in certain cases, such as when outliers are present. To address this issue, many statisticians have conducted research on bias reduction and have proposed variants of ML estimators based on their findings. Section \ref{sec:bias_reduction} will discuss the analogue of bias reduction for the proposed estimator.

As in the MD method used for the continuous response setup, we shall define a distance function that serves as the counterpart of the log-likelihood in the ML method. To this end, we define a CvM-type $L_{2}$ distance function
\ben\lel{eq:cL}
\cL(\B) = \sum_{j=1}^{J}\left[ \stk d_{kj}\Big\{\textrm{I}(y_{k}=1)-p_{k}(\B) \Big\}\right]^{2},
\een
where $\textrm{I}(\cdot)$ is an indicator function, $d_{kj}\in \R$, and $p_{k}(\B):=p(\bx_{k}'\B)$. Subsequently, define an MD estimator $\widehat{\B}$ as
\ben\lel{eq:opt}
\cL(\widehat{\B}) = \inf_{\B\in\R^{J}}\cL(\B).
\een
\begin{rem}\lel{rem:differentiability}
It is worth noting that $\textrm{I}(y_{k}=1)=y_{k}$ since $y_{k}$ is a binary variable. Hence, one can rewrite the distance function without the indicator function, leading to a simpler form; however, replacing the indicator function is not a viable option for a continuous response variable. In the continuous setup, distance functions used in the MD estimation method contain parameters inside the indicator function: e.g., $\textrm{I}(y_{k}-\bx_{k}'\B\le u)-H(u)$ where $u\in\mathbb{R}$, and $H$ is a df of some continuous random variable. Hence, the smoothness -- more precisely, the differentiability -- of the distance function with respect to (w.r.t) $\B$ cannot be affirmed, thereby rendering the optimization problem in (\ref{eq:opt}) rather challenging. On the contrary, the distance function in (\ref{eq:cL}) for binary regression equipped with the aforementioned $p$ (such as a logistic, normal, or Cauchy df) is now continuously differentiable w.r.t $\B$, which facilitates faster computation of the MD estimator than otherwise would be the case, thereby alleviating the computational burden of the MD estimation. Despite its desirable properties, the MD estimation has borne the brunt of criticism due to its computational cost; hence, its promising capability to deliver fast estimation will reinforce the merit of the MD methodology.
\end{rem}

Note that optimization problems with the CvM-type distance functions usually don’t have any closed-form solutions, and hence, finding an estimator tends to be somewhat daunting at first glance. Even after one somehow manages to procure the estimator, investigating its asymptotic properties will be a more challenging task. Addressing these issues, Koul (2002) introduced some conditions that are referred to as \textit{uniformly locally asymptotic quadratic} (ULAQ) conditions under which the distance function can be closely approximated by another quadratic distance function $\cQ$, and hence, the extent of difficulty in solving the complex optimization problem in (\ref{eq:opt}) can be pared down back to that in solving an optimization problem of the quadratic function. More importantly, moving pari passu with providing the estimator, the ULAQ conditions lays theoretical foundations upon which the asymptotic properties of the estimator can be derived. Thus, adopting the ULAQ conditions facilitates the applicability of the MD estimation to other problems and, as a result, expands the scope of the MD methodology. Before stating the ULAQ conditions, as in the literature on the MD estimation, the necessity of defining a neighborhood arises. For the true binary regression parameter $\B_{0}$, define its neighborhood as
\ben\lel{eq:neighbood}
\cN_{b}(\B_{0}):= \{\B\in \R^{p}:\|\bA^{-1}(\B-\B_{0})\|\le b\},
\een
where $\bA$ is a $J\times J$ symmetric, nonsingular matrix. More discussion of $\bA$ will appear in the following sections: e.g., see (a.1) in Section \ref{sec:ULAQ2}. Now we are ready to state the ULAQ conditions.
\begin{itemize}
\item[(\tbf{U.\,\addtocounter{Qcounter}{1}\theQcounter})] There exists a sequence of random vectors $\mbf{S}_{n}(\B_{0})\in \mathbb{R}^{J}$ and a sequence of $J\times J$ real matrices $\bW_{n}(\B_{0})$ such that for all $0<b<\iny$
\benn
\sup \left|\cL(\B)-\cL(\B_{0}) - (\B-\B_{0})'\mbf{S}_{n}(\B_{0}) - \frac{1}{2}(\B-\B_{0})'\bW_{n}(\B_{0})(\B-\B_{0})\right|=o_{p}(1),
\eenn
where the supremum is taken over $\cN_{b}(\B_{0})$.
\item[(\tbf{U.\,\addtocounter{Qcounter}{1}\theQcounter})]  For all $\vep>0$, there exists a $0<c_{\vep}<\iny$ such that
    \benn
    \mathbb{P}\left( |\cL(\B_{0})|\le c_{\vep} \right)\ge 1-\vep.
    \eenn
\item[(\tbf{U.\,\addtocounter{Qcounter}{1}\theQcounter})] For all $\vep>0$ and $0<c<\iny$, there exists a $0<b<\iny$ and $N_{b,\vep}$ (depending on $b$ and $\vep$) such that
    \benn
    \mathbb{P}\left(\inf |\cL(\B)| >c \right)\ge 1-\vep,
    \eenn
    where the infimum is taken over $\{\B\in \R^{p}:\|\bA^{-1}(\B-\B_{0})\|> b\}$.
\end{itemize}
\begin{rem}\lel{rem:ULAQ}
The first ULAQ condition implies the local, asymptotic quadraticity of $\cL$, that is, the distance function can be asymptotically approximated by a quadratic function (a sum of all other terms in the supremand) in the neighborhood of the true parameter: see (\ref{eq:cQ}). The other two conditions intuitively illustrate that the local, asymptotic quadraticity of the distance function is uniform: $\cL$ is stochastically bounded in the neighborhood, while it stochastically explodes to $\iny$ in a complement of the neighborhood, indicating that searching for the optimal solution to (\ref{eq:opt}) will amount to minimizing the quadratic function in $\cN_{b}(\B_{0})$.
\end{rem}

After ascertaining that the ULAQ conditions are met, one can obtain the asymptotic distribution of the MD estimator by using the following lemma, initially proposed by Koul (2002) and reproduced here: see Theorem 5.4.1 from Koul (op.\,\,cit).
\begin{lem}\lel{lem:asym_distr_Sn}
Assume that $\cL$ satisfies the ULAQ assumptions (\tbf{U.\,1})- (\tbf{U.\,\theQcounter}). Let $\cB_{n}:=\bA\bW_{n}(\B_{0})\bA$ where $\bA$ is as in (\ref{eq:neighbood}). Let $\wh{\B}$ denote the MD estimator that solves the optimization problem in (\ref{eq:opt}). Then,
\benn
\cB_{n}\bA^{-1} (\wh{\B}-\B_{0}) = -\bA\mbf{S}_{n}(\B_{0})+o_{p}(1).
\eenn
\end{lem}

\begin{rem}
Lemma \ref{lem:asym_distr_Sn} says finding the asymptotic distribution of the MD estimator $\wh{\B}$ amounts to specifying $\mbf{S}_{n}$, $\bW_{n}$, and the asymptotic distribution of $\mbf{S}_{n}$.
\end{rem}

\subsection{ULAQ conditions and asymptotic properties of the MD estimation}\lel{sec:ULAQ2}
\setcounter{Qcounter}{0}
This section focuses on two tasks at the top of the agenda: demonstrating that the ULAQ conditions for the distance function $\cL$ are satisfied, and identifying the asymptotic distribution of $\wh{\B}$ under these conditions. To achieve the goal, we introduce additional assumptions that are rooted in Koul (2002, pp. 173-175). This study requires fewer assumptions than Koul (2002), as some of the assumptions in Koul (op.\,\,cit.) are a priori satisfied; see, e.g., Remark \ref{rem:apriori}. Consider $n$ pairs of observations, $(y_{1}, \mbf{x}_{1}),...,(y_{n}, \mbf{x}_{n})$ where $\mbf{x}_{k}'\in \R^{J}$. Let $\bX$
be an $n\times J$ matrix the $k$th row vector of which is  $\mbf{x}_{k}$. Write a $n\times J$ matrix $\bD:=((d_{kj}))$, $1\le k\le n$, $1\le j\le J$, where $d_{kj}$'s are the real-valued weights used in (\ref{eq:cL}). Let $\mbf{d}_{k}'\in \R^{p}$, $1\le k\le n$ denote the $k$th row vector of $\bD$.
\begin{itemize}
\item[(\tbf{a.\addtocounter{Qcounter}{1}\theQcounter})] Let $\textbf{B}$ denote an $n\times n$ symmetric, positive definite matrix. Then, $(\bX'\textbf{B}\bX)^{-1}$ is nonsingular and $\bA:=(\bX'\textbf{B}\bX)^{-1/2}$ exists.
\item[(\tbf{a.\addtocounter{Qcounter}{1}\theQcounter})] For all $1\le j\le J$, $\sum_{k=1}^{n}d_{kj}^{2}=1$, and $\max_{1\le k\le n}d_{kj}=o(1)$.
\item[(\tbf{a.\addtocounter{Qcounter}{1}\theQcounter})] Let $\mbf{c}_{k}:=\bA \mbf{x}_{k}$. Then $\max_{1\le k\le n}\|\mbf{c}_{k}\|=o(1)$.
\item[(\tbf{a.\addtocounter{Qcounter}{1}\theQcounter})] For $1\le j\le J$, $\stk \|d_{kj}\mbf{c}_{k}\| = O(1)$.
\item[(\tbf{a.\addtocounter{Qcounter}{1}\theQcounter})] Let $f(s):=dp(s)/ds$. Then, $\|f\|_{\iny}<\iny$ and $\|f'\|_{\iny}<\iny$.
\item[(\tbf{a.\addtocounter{Qcounter}{1}\theQcounter})] For all $1\le k\le n$ and for all unit vectors $\mbf{e}\in \R^{J}$,  either $\mbf{d}_{k}'\mbf{e}\mbf{x}_{k}'\bA\mbf{e}\ge 0$ or $\mbf{d}_{k}'\mbf{e}\mbf{x}_{k}'\bA\mbf{e}\le 0$ holds true.
\end{itemize}
\begin{rem}
Note that (\tbf{a.5}) is met when $p$ is a logistic, normal, or Cauchy df.
\end{rem}
\begin{rem}
At first glance, assumption (\tbf{a.6}) is quite strong. However, it will be trivially met when $\bD=\bX\bA$, that is,
\benn
\mbf{d}_{k}'\mbf{e}\mbf{x}_{k}'\bA\mbf{e}=(\mbf{d}_{k}'\mbf{e})^{2}>0. \eenn
\end{rem}
Under (\tbf{a.1})-(\tbf{a.\theQcounter}), we shall prove that the distance function $\cL$ satisfies the ULAQ conditions that are introduced in the previous section and derive the asymptotic distribution of the MD estimator. Let $\W:=(\cW_{1},...,\cW_{p})'\in \R^{p}$ where
\benn
\cW_{j}(\B) = \sum_{k=1}^{n} d_{kj}\Big\{\textrm{I}(y_{k}=1)-p_{k}(\B) \Big\}.
\eenn
Observe that the distance function can be rewritten as
\benn
\cL(\B) = \sum_{j=1}^{J}\big[ \cW_{j}(\B)\big]^{2}=\|\W(\B)\|^{2}.
\eenn
Let $\mbf{q}_{k}(\B):=\partial p_{k}(\B)/\partial \B\in \R^{J}$. Subsequently, let
\benr
\mbf{S}_{n}(\B) &:=& -2\stj\stk  d_{kj}\cW_{j}(\B)\mbf{q}_{k}(\B),\lel{eq:Sn}\\
\bW_{n}(\B) &:=& 2\stj\sti\stk d_{kj}d_{ij}\mbf{q}_{k}(\B)\mbf{q}_{i}'(\B),\nonumber
\eenr
and
\ben\lel{eq:cQ}
\cQ(\B) := \cL(\B_{0}) + (\B-\B_{0})'\mbf{S}_{n}(\B_{0})+\frac{1}{2}(\B-\B_{0})'\bW_{n}(\B_{0})(\B-\B_{0}).
\een
Note that the verification of (\tbf{U.1}) amounts to showing for $0< b<\iny$
\benn
\sup_{\|\bA^{-1}(\B-\B_{0})\|\le b}|\cL(\B)-\cQ(\B)|=o_{p}(1).
\eenn
To prove the claim, we need the following lemma.
\ifnum \value{show}=1{
\noi
\begin{leftbar}
\noi
\tbf{Start} the proof of
\benn
\sup\stj\left[\stk d_{kj}\left\{ p_{k}(\B)-p_{k}(\B_{0}) -  (\B-\B_{0})'\mbf{q}_{k}(\B_{0}) \right\} \right]^2 = o(1).
\eenn
\noi
Refer to the \tecr{page 131} of Galaxy.
\end{leftbar}
}\fi
\begin{lem}\lel{lem:mu}
For $0< b<\iny$,
\ben\lel{eq:lem_mu}
\sup\stj\left[\stk d_{kj}\left\{ p_{k}(\B)-p_{k}(\B_{0}) -  (\B-\B_{0})'\mbf{q}_{k}(\B_{0}) \right\} \right]^2 = o(1),
\een
where the supremum is taken over $\|\bA^{-1}(\B-\B_{0})\|\le b$.
\end{lem}
\begin{proof}
Let $f_{k}(\B):=f(\mbf{x}_{k}'\B)$ where $f$ is as in (\tbf{a.5}). To conserve space, let $f_{k}^{0}$ and $f_{k}$ denote $f_{k}(\B_{0})$ and $f_{k}(\B)$, respectively, that is,
\benn
f_{k}^{0}=f(\mbf{x}_{k}'\B_{0}).
\eenn
Let $\mbf{u}:=\bA^{-1}(\B-\B_{0})\in \R^{J}$. Recall $\mbf{c}_{k}=\bA \mbf{x}_{k}\in \R^{J}$, $1\le k\le n$, and $\mbf{q}_{k}(\B)=f_{k}(\B)\mbf{x}_{k}'$. Thus,
\benn
p_{k}(\B)-p_{k}(\B_{0}) -  (\B-\B_{0})'\mbf{q}_{k}(\B_{0})
= \mbf{u}'\mbf{c}_{k}[f_{k}(\wt{\B})-f_{k}(\B)],
\eenn
where $\wt{\B}=c\B_{0}+(1-c)\B$ for some $c\in (0,1)$.
Note that
\benrr
|f_{k}(\wt{\B})-f_{k}(\B_{0})|
&\le& \|f'\|_{\iny} |\mbf{x}_{k}'(\B-\B_{0})|\\
&\le&\|f'\|_{\iny}\|\mbf{u}\|\|\mbf{c}_{k}\|,
\eenrr
where the first inequality follows from the mean value theorem and $|\mbf{x}_{k}'(\wt{\B}-\B)|\le |\mbf{x}_{k}'(\B-\B_{0})|$. Therefore,
\benrr
\textrm{the supremand of LHS of (\ref{eq:lem_mu})}
&\le&  \stj \left[ \stk |d_{kj}\mbf{u}'\mbf{c}_{k}|\cdot |f_{k}(\wt{\B})-f_{k}(\B_{0})| \right]^{2} \\
&\le& b^4 \|f'\|_{\iny}^2 \left(\max_{1\le k\le n}\|\mbf{c}_{k}\|\right)^2\stj\left(\stk \|d_{kj}\mbf{c}_{k}\|\right)^2\longrightarrow 0,
\eenrr
where (\tbf{a.3}), (\tbf{a.4}) and (\tbf{a.5}) imply the convergence of zero, thereby completing the proof of the claim.
\end{proof}
\noi
The next theorem demonstrates the first ULAQ condition is indeed satisfied.
\begin{thm}\lel{thm:L_Q}
Assume (\tbf{a.1})-(\tbf{a.\theQcounter}). Then, the distance function $\cL$ in (\ref{eq:cL}) satisfies \textbf{(U.1)}, that is, for any $0<b<\iny$,
\benn
   \mathbb{E}\Big(\sup|\cL(\B)- \cQ(\B)|\Big)=o(1),
\eenn
where the supremum is taken over $\{\B\in\R^{p}:\|\bA^{-1}(\B-\B_{0})\|\leq b\}$.
\end{thm}
\noi

\begin{proof} %The proof will be similar to one for Theorem 1.1 in Kim (2024), albeit more complicating.
Consider $\mbf{u}:=\bA^{-1}(\B-\B_{0})$ with $\|\mbf{u}\|\le b<\iny$. Define $\mbf{R}_{j}(\B):=\stk d_{kj}f_{k}(\B)\mbf{x}_{k}$, and let $\bR(\B):=[\mbf{R}_{1}(\B),...,\mbf{R}_{J}(\B)]$ denote a $J\times J$ matrix whose $j$th column is $\mbf{R}_{j}(\B)$. Note that $\cL$ and $\cQ$ can be rewritten in the following quadratic forms
\benrr
\cL(\B) &=& \stj\Bigg[\Big\{\cW_{j}(\B_{0}) - (\B-\B_{0})'\mbf{R}_{j}(\B_{0})   \Big\} - \left. \stk d_{kj}\Big\{ p_{k}(\B)-p_{k}(\B_{0}) - (\B-\B_{0})'\mbf{R}_{j}(\B_{0}) \Big\} \right]^{2},
\eenrr
and
\benrr
\cQ(\B) &=& \stj\left[\cW_{j}( \B_{0})-(\B-\B_{0})'\mbf{R}_{j}(\B_{0}) \right]^{2}= \|\W(\B_{0})-(\B-\B_{0})'\bR\|^2.
\eenrr
\ifnum \value{show}=1{
Note that
\benn
\mathbb{E} \|\W(\B_{0})\|^{2} = \sum_{j=1}^{J}\stk d_{kj}^2p_{k}(\B)(1-p_{k}(\B))\le \sum_{j=1}^{J} \stk \le J. d_{kj}^2
\eenn
}\fi
$\mathbb{E} \|\W(\B_{0})\|^{2}<\iny$ readily follows from (\tbf{a.2}), which, in turn, implies
\ben\lel{eq:Op1}
\|\W(\B_{0})\|^{2} = O_{p}(1).
\een
Next, observe that
\benr
(\B-\B_{0})'\bR\bR'(\B-\B_{0})
&=&  \stj |\mbf{u}'\bA\mbf{R}_{j}(\B_{0})|^2  \lel{eq:O1}\\
&\le& b^2 \|f\|_{\iny}^{2} \stj\left(\stk \| d_{kj}\mbf{c}_{k}\|\right)^2 <\iny, \nonumber
\eenr
where the finiteness immediately follows from (\tbf{a.4}) and (\tbf{a.5}). Then, in view of Lemma \ref{lem:mu}, (\ref{eq:Op1}) and (\ref{eq:O1}), expanding the quadratic of $\cL$ and applying the C-S inequality to the cross product term will prove the claim, thereby completing the proof of the theorem.
\end{proof}

\begin{rem}\lel{rem:apriori} In the case of a continuous response variable $y$, Koul (2002, pp.\,173-175) assumed  (\ref{eq:Op1}) and showed that several pairs of integrating measure and distribution satisfy the assumption, which is a priori met in this study.
\end{rem}
\noi
The following lemma verifies that the remaining ULAQ conditions are satisfied.
\begin{lem}\lel{lem:asym_bound}
Suppose the assumptions in Theorem \ref{thm:L_Q} hold. Then, $\cL$ also satisfies (\tbf{U.2}) and (\tbf{U.3}).
\end{lem}
\begin{proof}
Let $\mbf{u}\in \mathbb{R}^{J}$. Define $\mbf{V}(\mbf{u}):=(V_{1}(\mbf{u}),...,V_{J}(\mbf{u}))'\in \R^{J}$ and $\widehat{\mbf{V}}(\mbf{u}):=(\widehat{V}_{1}(\mbf{u}),...,\widehat{V}_{J}(\mbf{u}))'\in \R^{J}$ where for $1\le j\le J$,
\benn
V_{j}(\mbf{u}):= \cW_{j}(\B_{0}+\bA\mbf{u})\textrm{ and }\widehat{V}_{j}(\mbf{u}):= \big\{\cW_{j}(\B_{0})+\mbf{u}'\bA \mbf{R}_{j}(\B_{0})\big\}.
\eenn
As in Lemma 5.5.4 in Koul (2002), we shall show that $\mbf{e}'\mbf{V}(\mbf{u})$ is monotone in $\|\mbf{u}\|$ together with the assumption (\tbf{a.6}). Let $\mbf{u}=r\mbf{e}$ with $\|\mbf{e}\|=1$. Note that
\benn
\mbf{e}'\mbf{V}(\mbf{u}) = \stk (\mbf{e}'\mbf{d}_{k})\left[\textrm{I}(y_{k}=1)-\int_{0}^{x_{k}'\beta_{0}+ re'\bA \mbf{x}_{k}}f(s) ds\right],
\eenn
where $\mbf{e}'\mbf{V}(\mbf{u})$ decreases as $r$ increases
if $\mbf{d}_{k}'\mbf{e}\mbf{x}_{k}'\bA\mbf{e}\ge 0$ for all $k$, while the opposite holds if the inequality changes. Recall $\G_{n}(\B)$ and let $k_{n}(\mbf{e}):=\mbf{e}'\G_{n}\mbf{e}$. Then, the rest of  the proof will be the same as the proof of Lemma 5.5.4 in Koul (2002) and is omitted here.
\end{proof}
\ifnum \value{show}=1{
\noi
\begin{leftbar}
\noi
\tbf{Mar/23/2025} Start the proof of
\benn
\bA\mbf{S}_{n}(\B_{0})\Rightarrow_{\cD}N(0, \wt{\S}).
\eenn
\noi
Refer to the \tecr{page 134} of Galaxy (Mar 19 2025).
\end{leftbar}
}\fi
Ascertaining that the ULAQ conditions are all met, we proceed to identify the asymptotic distribution of $\wh{\B}$. Recall $\bW_{n}(\B)$ in (\ref{eq:cQ}). Define an $n\times n$ diagonal matrix $\L(\B):=\textrm{diag}[f_{1}(\B),...,f_{n}(\B)]$ and a $J\times n$ matrix $\G_{n}(\B):=\bA\bX'\L(\B)\bD$. Then, it is not difficult to see that $\bR=\bA^{-1}\G_{n}$, and hence, one readily has
\benrr
\bW_{n}(\B)&=&\stj \mbf{R}_{j}\mbf{R}_{j}',\\
&=& \mathbf{A}^{-1}\wt{\G}_{n}(\B)\mathbf{A}^{-1},
\eenrr
where $\wt{\G}_{n}(\B):=\G_{n}(\B)\G_{n}(\B)'$. Note that
$\wt{\G}_{n}(\B)$ will play the role of $\cB_{n}$ in Lemma \ref{lem:asym_distr_Sn}. Let $\si_{ij}:=\stk p_{k}(1-p_{k})d_{ki}d_{kj}$ and $\S_{n}$ denote a $J\times J$ matrix whose $(i,j)$th entry is $\si_{ij}$. The asymptotic distribution of $\wh{\B}$ will directly follow from the next lemma.
\begin{lem}\lel{lem:Asym_Convergence}
Assume that $\S_{n}$ is positive definite and
\ben\lel{eq:sig}
\lim_{n\ra\iny}\G_{n}(\B_{0})\S_{n}\G_{n}'(\B_{0}) = \wt{\S},
\een
where $\wt{\S}$ possibly depends on $\B_{0}$. Then,
\benn
\bA\mbf{S}_{n}(\B_{0})\Rightarrow_{\cD}N(0, 4\wt{\S}).
\eenn
\end{lem}
\begin{proof}
Recall $\mbf{S}_{n}$ and rewrite $\mbf{S}_{n}(\B_{0})=-2\bR(\B_{0})\W(\B_{0})$. Observe that
\benn
\bA \mbf{S}_{n}(\B_{0})=-2\bA\bX'\L(\B)\bD\W=-2\G_{n}(\B_{0})\W.
\eenn
Thus, the proof of the lemma amounts to showing asymptotic convergence of $\W$ to a normal distribution. Let $\eta_{k}(\B):=I(Y_{k}=1)-p_{k}(\B)$. Note that $\mathbb{E}[\eta_{k}]=0$ and $\mathbb{E}[\eta_{k}^{2}]=p_{k}(1-p_{k})$. For $\mbf{a}:=(a_{1},...,a_{J})\in \R^{J}$, one has $\mbf{a}'\W=\stk t_{k}$ where
\benn
t_{k}:=\eta_{k}\sum_{j=1}^{J} a_{j}d_{kj}=\eta_{k}\tilde{a}_{k}, \quad (say).
\eenn
Note that from $|\eta_{k}|\le 1$, we have
\ben\lel{eq:tk}
|t_{k}|\le |\tilde{a}_{k}|,
\een
and \tbf{(a.2)} implies that
\ben\lel{eq:ak}
\max|\tilde{a}_{k}|=o(1).
\een
Let $\tau_{n}^{2}:=\stk \mathbb{E}[t_{k}^{2}]$. It is not difficult to show that
\ben\lel{eq:tau}
\tau_{n}^{2} = \stk \tilde{a}_{k}^{2}p_{k}(1-p_{k}).
\een
We shall show that the Lindeberg-Feller (LF) condition for the convergence of $\mbf{a}'\W$ is met, that is, for all $\vep>0$
\benrr
\tau_{n}^{-2}\stk \mathbb{E}\left[t_{k}^{2}:|t_{k}|>\vep \tau_{n}\right]&\le& \tau_{n}^{-2} \stk \tilde{a}_{k}^{2}\mathbb{P}(|t_{k}|>\vep \tau_{n})\\
&\leq& \tau_{n}^{-4}\vep^{-2} \stk \tilde{a}_{k}^{2}\mathbb{E}(t_{k}^{2}) \\
&\le&\vep^{-2}\tau_{n}^{-2} \max_{1\le k\le n}\tilde{a}_{k}^{2} \longrightarrow 0,
\eenrr
where the first, second, and last inequalities follow from (\ref{eq:tk}), the Chevyshev inequality and (\ref{eq:tau}), respectively, while the assumption (\tbf{a.2}) and (\ref{eq:ak}) readily imply the convergence to 0. Thus, one has
\benn
\tau_{n}^{-1}\stk t_{k}\Rightarrow N(0,1).
\eenn
Note that $\tau_{n}^{2}=\sum_{i}\sum_{j}a_{i}a_{j}\si_{ij}=\mbf{a}'\S_{n}\mbf{a}$, and hence, (\ref{eq:sig}) and the application of Cramer-Wold device yield
\benn
\S_{n}^{-1/2}\W\Rightarrow_{\cD}N(\mathbf{0}_{J\times 1}, \mathbf{I}_{J\times J}),
\eenn
thereby completing the proof of the lemma.
\end{proof}

\noi
We conclude this section by stating the main result of this study: the asymptotic distribution of the MD estimator. Recall $\wt{\G}_{n}(\B)$.
\begin{thm}\lel{thm:asymp}
Suppose the assumptions in Theorem \ref{thm:L_Q} and Lemma \ref{lem:Asym_Convergence} hold. Then
the MD estimator $\wh{\B}$ asymptotically follows the normal distribution, that is,
\benn
\wt{\G}_{n}(\B)\bA^{-1}(\wh{\B}-\B_{0}) \Rightarrow_{\cD} N(\mathbf{0}, 4\wt{\S}),
\eenn
where $\wt{\S}$ is as in Lemma \ref{lem:asym_distr_Sn}.

\end{thm}

\begin{rem}
Note that the matrix $\wt{\S}$ is a function of the unknown $\B_{0}$, and hence, any statistical inferences about $\B_{0}$ will heavily depend on the MD estimation: for example, the standard error of the estimator should be estimated by using $\wh{\B}$. Due to the consistency of $\wh{\B}$ and the fact that all entries of $\wt{\S}$ are continuous in $\wh{\B}$, the convergence of $\wt{\S}(\wh{\B})$ to $\wt{\S}(\B_{0})$ in probability can be easily established. However, there still exists a chance that any procedures that rely heavily on the estimate lead to incorrect statistical inference: see, e.g., Mansournia (2018).
\end{rem}

\begin{proof} Theorem \ref{thm:L_Q} and Lemma \ref{lem:asym_bound} imply the distance function $\cL$ satisfies the ULAQ conditions. Consequently, Lemmas \ref{lem:asym_distr_Sn} and \ref{lem:Asym_Convergence} yield the desired result, thereby completing the proof of the theorem.
\end{proof}

\ifnum \value{show}=1{
\noi
\begin{leftbar}
\noi
\tbf{April/16/2025} Asymptotic variance of the MDE.
Refer to the \tecr{page 28} of ``Robust part" (April 16, 2025).
\end{leftbar}
}\fi

\begin{rem}\lel{rem:best}
Let $\mathbf{P}_{n}:=diag[p_{1}(\B)(1-p_{1}(\B)),...,p_{n}(\B)(1-p_{n}(\B))]$. Note that $\S_{n}=\bD'\mathbf{P}_{n}\bD$. Let $AVar(\wh{\B})$ denote the asymptotic variance of $\wh{\B}$. Then, $AVar(\wh{\B})$ can be written as
\benn
\bA(\wt{\G}_{n}^{-1})'\wt{\S}\wt{\G}_{n}^{-1}\bA=\bA(\G_{n}')^{-1}\bD' \mathbf{P}_{n}\bD \G_{n}^{-1}\bA.
\eenn
Recall $\L(\B):=\textrm{diag}[f_{1}(\B),...,f_{n}(\B)]$ and choose $\bA$ such that $\bA^{2}=(\bX'\L \mathbf{P}_{n}^{-1}\L\bX)^{-1}$. Observe that $AVar(\wh{\B})$ can be rewritten as
\benn
AVar(\wh{\B}) = \bA^{2} = (\bX'\L \mathbf{P}_{n}^{-1}\L\bX)^{-1}.
\eenn
\end{rem}

\begin{rem}
Consider the logistic $p(s)=1/(1+e^{-s})$. For $\mbf{u}:=(u_{1},...,u_{n})'\in \mathbb{R}^{n}$, let $\mathbf{e}^{\mbf{u}}:=diag[e^{u_{1}},...,e^{u_{n}}]$. Note that the $i$th entry of the diagonal matrix $\L \mathbf{P}_{n}^{-1}\L$ is $e^{\bx_{i}'\B}$, and hence,
\benn
AVar(\wh{\B})=(\bX'\mathbf{e}^{\bX'\B}\bX)^{-1}.
\eenn
\end{rem}

\begin{rem}
Observe that both $\mathbf{P}_{n}$ and $\L$ are matrices whose entries are functions of $\B$. Thus, finding $AVar(\wh{\B})$ requires the information of the true $\B$ that is unknown. Replacing $\B$ with the MDE, which is a consistent estimator for $\B$, we can obtain a consistent estimator for $AVar(\wh{\B})$.
\end{rem}

\subsection{Robustness}\lel{Sec:robustness}

\ifnum \value{show}=1{
\noi
\begin{leftbar}
\noi
\tbf{April/13/2025} Start the Efficiency and Robustness
Refer to the \tecr{page 36} of Galaxy (April 17 2025).
\end{leftbar}
}\fi

Recall $n$ pairs of the observations $\{(y_{i}, \bx_{i})'\in \mathbb{R}^{(J+1)}:1\le i\le n\}$ and the log-likelihood of the ML estimation in (\ref{eq:log_likelihood}), where $\mathbb{P}(y_{i}=1)=p(\bx_{i}'\B)$ with $f(u)=dp(u)/du$. As done in Copas (1988), the ML estimator can be obtained by solving
\ben\lel{eq:sum_psi}
\sti \Psi(y_{i}, \bx_{i}, p)=0,
\een
where
\benn
\Psi(y_{i}, \bx_{i}, p):=\frac{f(\bx_{i}'\B)[y_{i}-p(\bx_{i}'\B)]}{p(\bx_{i}'\B)[1-p(\bx_{i}'\B)]}\bx_{i}.
\eenn
Let $\Delta_{L}(y_{i},\bx_{i},p)$ denote the influence function of the ML estimation at $(y_{i}, \bx_{i})$. Then, being analogous to (2.3.5) of Hampel et al.\,\,(1986, p.101), the influence function can be written as
\benn
\Delta_{L}(y_{i},\bx_{i},p) = \left\| \frac{(y_{i}-p(\bx_{i}'\B))}{f(\bx_{i}'\B)}(\bx_{i}\bx_{i}')^{-1}\bx_{i}\right\|,
\eenn
where $(\bx_{i}\bx_{i}')^{-1}$ is a $J\times J$ matrix, and $\|\cdot\|$ is the Euclidean norm.
Consider $J=1$. Then,
\benn
\Delta_{L}(y_{i},x_{i},p)=\frac{|y_{i}-p(\beta x_{i})|}{|x_{i}f(\beta x_{i})|}.
\eenn
Note that when $f$ is the logistic, normal, or Cauchy density function, $uf(cu)$ will converge to 0 for $c\neq0$ as $u$ goes to $\iny$. Hence, the influence function is unbounded, implying that the ML estimator is vulnerable to the presence of outliers.

Next, consider the MD estimator. Recall $\cW_{j}(\B)$, $1\le j\le n$, and let $\cW_{i}^{*}(\B):=\sum_{j=1}^{J}d_{ij}\cW_{j}(\B)$ where $d_{ij}\in \R$ is the $(i,j)$th entry of $\bD$. Then, the the MD estimator can be obtained from solving the equation in (\ref{eq:sum_psi}) where
\benn
\Psi(y_{i}, \bx_{i}, p)= f(\bx_{i}'\B)\cW_{i}^{*}(\B)\bx_{i}.
\eenn
Consequently, the influence function of the MD estimator -- denoted by $\Delta_{D}$ -- at $(y_{i},\bx_{i})$ will be
\benn
\Delta_{D}(y_{i},\bx_{i},p)= \left\|\cW_{i}^{*}(\B)\left(\stk d_{ki}^{*} f(\bx_{k}'\B)\bx_{k}\bx_{i}' \right)^{-1}\bx_{i}\right\|,
\eenn
where $d_{ki}^{*}$ is the $(k,i)$th entry of the matrix $\bD^{*}=\bD\bD'$: see Hampel (1986) and Kim (2025) for more details of how to obtain the influence function of the MD estimator. In the case of $J=1$, the influence function of the MD estimation will be simplified as follows:
\benn
\Delta_{D}(y_{i},x_{i},p)=\frac{\left|\cW_{i}^{*}(\B)\right|}{\sum_{k\neq i}d_{ki}^{*} f(\beta x_{k})x_{k}+d_{ii}^{*} f(\beta x_{i})x_{i}}.
\eenn
As mentioned earlier, extremely large $x_{i}$ will cause $f(\beta x_{i})x_{i}$ to be approximately zero; however, the denominator of $\Delta_{D}$ is apart  from 0 as much as $|\sum_{k\neq i}d_{ki}^{*} f(\beta x_{k})x_{k}|$, and hence, $\Delta_{D}$ remains finite, even in the presence of the extremely large $x_{i}$. Unlike the influence function of the ML estimator, that of the MD estimator is bounded. Hence, the impact of any outliers will be limited, thereby leading the MD estimator to be robust to outliers.

\subsection{Bias reduction}\lel{sec:bias_reduction}
As will be shown in the simulation studies, the MD method reports a relatively larger bias but a much smaller standard error (SE) than other competing methods. Fortunately, its root mean squared error (RMSE) turns out to be much better -- due to the smaller SE -- than other methods. Thus, the MD method will consolidate its superiority over others if its bias issue is addressed. For this reason, we will rigorously investigate whether bias reduction of the MD estimator is feasible.

The statistical literature is replete with research articles discussing bias reduction various some estimators. One of the popular compendia for reference includes Cox and Snell (1968). They expanded the score function of the likelihood to the second order, used an approximation, and demonstrated that the resulting estimator has bias of order $O(n^{-1})$, where $n$ denotes the sample size. Copas (1988) applied the analogous method -- the expansion and approximation of the score function to the second order -- to the binary regression model and showed that the bias remains of order $O(n^{-1})$. To obtain a less biased ML estimator, Firth (1993) proposed the modified score function that is downward-shifted from the original score function by bias multiplied by the Fisher information; he also demonstrated that in exponential families, the optimally penalized likelihood turns out to be the original likelihood factored by Jeffreys invariant prior, and the $O(n^{-1})$ term of the bias of the ML estimator is successfully removed. Kosmidis and Firth (2009) extended the application of the bias correction to a broader class of generalized nonlinear models. Focusing on the binary regression model alone, Kosmidis and Firth (2021) showed that using the Jeffreys prior as a penalty function reduces the asymptotic bias of ML estimators across various link functions, including logistic, probit, and log-log.

Analogously, we will expand the first derivative of $\cL$, which is a counterpart of the score function of the ML estimation, to the second order and explicitly express the bias of the MD estimator in an equation. Let $\partial \cL(\B)/\partial \B$ and $\partial^{i} \cL(\B)/\partial \B^{i}$ denote the first and $i$th order derivatives of $\cL$ with respect to $\B$, respectively. Using the Taylor expansion, we have
\ben\lel{eq:taylor_l}
\left. \frac{\partial\cL(\B)}{\partial \B}  \right|_{\B=\wh{\B}}
\approx \left. \frac{\partial\cL(\B)}{\partial \B}  \right|_{\B=\B_{0}}+(\wh{\B}-\B_{0})'\left. \frac{\partial^{2}\cL(\B)}{\partial \B^{2}}  \right|_{\B=\B_{0}}+\frac{1}{2}(\wh{\B}-\B_{0})'
\left. \frac{\partial^{3}\cL(\B)}{\partial \B^{3}}  \right|_{\B=\B_{0}}
(\wh{\B}-\B_{0}).
\een
Using (\ref{eq:opt}) together with the fact that the expectation of the first term in the right-hand side (RHS) of the above equation is actually zero, we shall derive the approximate bias from the above equation after taking the expectation on both sides thereof. Recall $\cW_{k}^{*}(\B)=\sum_{j=1}^{J}d_{kj}\cW_{j}(\B)$ from Section \ref{Sec:robustness}. Note that
\benrr
\frac{\partial\cL(\B)}{\partial \B}
&=& -2\stk \cW_{k}^{*}(\B)\mbf{q}_{k}(\B),\\
&=& -2\stk \mbf{U}^{k}(\B),\qquad (say).
\eenrr
For $1\le j\le J$, let $U_{j}^{k}$ denote the $j$th entry of $\mbf{U}^{k}\in\mathbb{R}^{J}$. For $1\le r\le J$, the partial derivative of $\cW_{k}^{*}$ with respect to $\beta_{r}$ is
\benn
\frac{\partial \cW_{k}^{*}}{\partial \beta_{r}}=-\sum_{l=1}^{n}d_{kl}^{*}f_{l}(\B)x_{lr},
\eenn
where $d_{kl}^{*}=\sum_{j=1}^{J}d_{kj}d_{lj}$, and hence, we have
\ben\lel{eq:first_derivative}
\frac{\partial U_{j}^{k}}{\partial \beta_{r}}=-f_{k}(\B)x_{kj}\sum_{l=1}^{n}d_{kl}^{*}f_{l}(\B)x_{lr}
+ \cW_{k}^{*}(\B)f_{k}'(\B)x_{kr}x_{kj}.
\een
In what follows, we write $f_{k}'(\B)$ and $f_{k}''(\B)$ as $f_{k}'$ and $f_{k}''$, respectively, to conserve space. As mentioned earlier, taking expectations on both sides of (\ref{eq:taylor_l}) will yield
%\ben\lel{eq:expansion_Uk}
%\mbf{U}^{k}(\wh{\B})\approx \mbf{U}^{k}(\B_{0})+(\wh{\B}-\B_{0})'
%\left. \frac{\partial \mbf{U}^{k}(\B)}{\partial \B}  \right|_{\B=\B_{0}}
%+\frac{1}{2}(\wh{\B}-\B)'
%\left. \frac{\partial^{2} \mbf{U}^{k}(\B)}{\partial \B^{2}}  \right|_{\B=\B_{0}}(\wh{\B}-\B).
%\een
\ben\lel{eq:expansion_Uk}
\stk\mathbb{E}\left[(\wh{\B}-\B_{0})'
\left. \frac{\partial \mbf{U}^{k}(\B)}{\partial \B}  \right|_{\B=\B_{0}}\right]
+\frac{1}{2}\stk\mathbb{E}\left[(\wh{\B}-\B_{0})'
\left. \frac{\partial^{2} \mbf{U}^{k}(\B)}{\partial \B^{2}}  \right|_{\B=\B_{0}}(\wh{\B}-\B_{0})\right]=0.
\een
Note that the $r$th entry of $\mathbb{E}\left[(\wh{\B}-\B)' (\partial \mbf{U}^{k}/\partial \B)  \right]$ will be
\benrr
\mathbb{E}\left[\sum_{i=1}^{J}(\hat{\beta}_{i}-\beta_{i})\frac{\partial U_{i}^{k}}{\partial \beta_{r}} \right]&=&-\sum_{i=1}^{J}\left[ f_{k}(\B)x_{ki}\stl d_{kl}^{*} f_{l}(\B)x_{lr}
\right]\mathbb{E}(\hat{\beta}_{i}-\beta_{i})+\sum_{i=1}^{J}\dot{f}_{k}(\B)x_{kr}x_{ki}\mathbb{E}\left[(\hat{\beta}_{i}-\beta_{i})\cW_{k}^{*}\right],\\
&=&-\sum_{i=1}^{J}c_{r1}^{ki}\mathbb{E}(\hat{\beta}_{i}-\beta_{i})+\sum_{i=1}^{J}c_{r2}^{ki}\mathbb{E}\left[(\hat{\beta}_{i}-\beta_{i})\cW_{k}^{*}\right],
\eenrr
where $c_{r1}^{ki}= f_{k}(\B)x_{ki}\stl d_{kl}^{*} f_{l}(\B)x_{lr}$ and $c_{r2}^{ki}= \dot{f}_{k}(\B)x_{kr}x_{ki}$.

Next, recall that
\ben\lel{eq:betahat}
(\hat{\B}-\B)=2\stk \bW_{n}^{-1}\mbf{U}^{k}(\B).
\een
Let $\mbf{w}_{i}'$ denote the $i$th row vector of $\bW_{n}^{-1}$, and write the $i$th entry of $(\hat{\B}-\B)$ as
\ben\lel{eq:betahat_i}
(\hat{\beta}_{i}-\beta_{i}) = 2\stk\cW_{k}^{*}\mbf{w}_{i}'\mbf{q}_{k}.
\een
A direct, albeit complex, calculation shows that
\benrr
\mathbb{E}\left[(\hat{\beta}_{i}-\beta_{i})\cW_{k}^{*}\right]&=&\sum_{j=1}^{J}d_{kj}Cov\left[(\hat{\beta}_{i}-\beta_{i}), \cW_{j}\right],\\
&=& 2\sum_{j=1}^{J}d_{kj} \stt \mbf{w}_{i}'\mbf{q}_{t} \sum_{h=1}^{J}d_{th} Cov(\cW_{h}, \cW_{j}),\\
&=& 2\sum_{j=1}^{J}d_{kj} \stt \mbf{w}_{i}'\mbf{q}_{t} \sum_{h=1}^{J}d_{th} \stl d_{lh}d_{lj}p_{l}(1-p_{l}),\\
&=& 2c_{3}^{ki}, \qquad (say),
\eenrr
where the second line follows from (\ref{eq:betahat_i}) and the definition of $\cW_{k}^{*}$.

Next, we proceed to check the last term on the left-hand side (LHS) of (\ref{eq:expansion_Uk}). Define $\Delta \mathbf{U}_{r}^{k}:=
\frac{\partial}{\partial \beta_{r}}\left(\frac{\partial U^{k}}{\partial \B}\right)\in \mathbb{R}^{J\times J}$.
Recall $\mbf{w}_{i}'$ and $\mbf{w}_{j}'$ that are the $i$th and $j$th row vectors of $\mathbf{W}_{n}^{-1}$, respectively. Using (\ref{eq:betahat}), the $r$th entry of the last term is, then, $2\sth\stl \left[(\mbf{U}^{h})'\mathbf{W}_{n}^{-1}\Delta\mathbf{U}_{r}^{k}\mathbf{W}_{n}^{-1}\mbf{U}^{l}\right]$. Let $a_{r}^{khl}$ denote the summand; it is not difficult to see that
\benn
a_{r}^{khl}=\sum_{i=1}^{J}\sum_{j=1}^{J} U_{j}^{h}U_{i}^{l}\sum_{s=1}^{J}\sum_{t=1}^{J}w_{js}w_{it}\frac{\partial^2 U_{s}^{k}}{\partial \beta_{r}\partial \beta_{t}},
\eenn
where $w_{it}$ and $w_{js}$ are the $t$th and $s$th entries of $\mbf{w}_{i}'$ and $\mbf{w}_{j}'$, respectively. Using (\ref{eq:first_derivative}) again and taking an expectation, we have
\benrr
\mathbb{E}(a_{r}^{khl})&=& \sum_{i=1}^{J}\sum_{j=1}^{J} \sum_{s=1}^{J}\sum_{t=1}^{J} w_{js}w_{it}
Cov\left( U_{j}^{h}, U_{i}^{l}\frac{\partial^2 U_{s}^{k}}{\partial \beta_{r}\partial \beta_{t}} \right).
\eenrr
Note that
\benn
Cov\left( U_{j}^{h}, U_{i}^{l}\frac{\partial^2 U_{ks}}{\partial \beta_{r}\partial \beta_{t}} \right) = \zeta_{1}Cov(U_{j}^{h}, U_{i}^{l})+\zeta_{2}Cov(U_{j}^{h}, U_{i}^{l}\cW_{k}^{*}),
\eenn
where
\benrr
\zeta_{1}&:=& -x_{ks}\sum_{m=1}^{n}d_{km}^{*}x_{mr}(f_{k}'x_{kt}f_{m}+f_{k}f_{m}'x_{mt}) -f_{k}'x_{kr}x_{ks}\sum_{m=1}^{n} d_{km}^{*}f_{k}x_{mt},\\
\zeta_{2}&:=&f_{k}''x_{kt}x_{kr}x_{ks}.
\eenrr
Direct calculations show that
\benrr
Cov(U_{j}^{h}, U_{i}^{l})&=&f_{h}f_{l}x_{hj}x_{li}\sum_{s=1}^{J}\sum_{t=1}^{J} d_{hs}d_{lt}\sum_{m=1}^{n}d_{ms}d_{mt}p_{m}(1-p_{m}),\\
Cov(U_{j}^{h}, U_{i}^{l}\cW_{k}^{*})&=& f_{h}f_{l}x_{hj}x_{li}\sum_{j=1}^{J}d_{hj}\sum_{s=1}^{J}\sum_{t=1}^{J}d_{ls}d_{kt}\kappa_{st}^{j},
\eenrr
where
\benn
\kappa_{st}^{j} = \sum_{m=1}^{n}d_{mj}d_{ms}d_{mt}(p_{m}-3p_{m}^{2}+2p_{m}^{3}).
\eenn
Let $\mbf{b}:=\mathbb{E}(\wh{\B}-\B)'\in \mathbb{R}^{J}$, which is a vector of the biases of the MD estimator. Let $b_{i}$ denote the $i$th entry of $\mbf{b}$. Finally, we obtain $\mbf{b}$ by solving
\benn
-\stk \sum_{i=1}^{J}c_{r1}^{ki}b_{i}+2\stk\sum_{i=1}^{J}c_{r2}^{ki}c_{3}^{ki}=\stk\sum_{h=1}^{n}\sum_{l=1}^{n} \mathbb{E}(a_{r}^{khl}).
\eenn
Note that the above equation can be written in matrix form, that is,
\benn
-\mathbf{C}\mbf{b}+2\mbf{d}=\frac{1}{2}\mbf{a},
\eenn
where $\mathbf{C}$ is a $J\times J$ matrix, and $\mbf{a}$ and $\mbf{d}$ are $J\times 1$ vectors. More specifically, for $1\le r, q\le J$, let $C_{rq}$, $a_{r}$, and $d_{r}$ denote the $(r,q)$th and $r$th entries of $\mathbf{C}$, $\mbf{a}$, and $\mbf{d}$, respectively. Then,
\benn
C_{rq} = \stk c_{r1}^{kq},\,\,
d_{r} = \stk\sum_{i=1}^{J}c_{r2}^{ki}c_{3}^{ki},\,\,
a_{r} = \stk \stl\sth \mathbb{E}(a_{r}^{khl}).
\eenn
It is not difficult to see that $\mathbf{C}=\bX'\L\bD^{*}\L\bX$ where $\bD^{*}:=\bD\bD'$. Next, let $\wt{\mathbf{W}}:=\mathbf{I}_{n\times n}\bigotimes(\mathbf{W}_{n}^{-1}\bX'\L)$ where $\mathbf{I}_{n\times n}$ is an $n\times n$ identity matrix, and $\bigotimes$ denotes the Kronecker product. Thus, $\wt{\mathbf{W}}$ will be an $nJ\times n^2$ matrix and can be partitioned into an $n\times n$ block so that all its diagonal blocks will be $\mathbf{W}_{n}^{-1}\bX'\L$, while all other blocks will be $n\times J$ zero matrices. Let $\mbf{P}_{k}$ denote the $k$th column vector of $\bD^{*}\mathbf{P}_{n}\bD^{*}$. Subsequently, define an $n^2\times 1$ vector $\wh{\mbf{P}}$ by stacking $\mbf{P}_{k}$, $k=1,2,...,n$. Finally, define a $J\times nJ$ matrix $\mathbf{M}$ whose $r$th row vector can be partitioned into $n$ pieces of $1\times J$ subvectors, where the $k$th subvector is $(c_{r2}^{k1},...,c_{r2}^{kn})$. Then $\mbf{d}$ can be rewritten as
\benn
\mbf{d}=\mathbf{M}\wt{\mathbf{W}}\wh{\mbf{P}}.
\eenn
Note that $\mbf{a}$ can also be expressed as a product of matrices and a single vector, as in $\mbf{d}$, even though the expression of $\mbf{a}$ will be a bit more complicated than that of $\mbf{d}$, and hence, it is not included here.

\section{Empirical studies}\lel{Sec:emp_studies}
When datasets are required for the simulation study in this section, we will generate $(y_{i},\mbf{x}_{i})'\in \mathbb{R}^{J+1}$, $1\le i\le n$ as follows. First, we generate $\mbf{x}_{i}=(x_{i1},...,x_{iJ})'\in \mathbb{R}^{J}$, $1\le i\le n$, where $x_{ij}$, $1\le j\le J$, are uniform random numbers between 0 and 3. Next, using a logistic $p(s)=1/(1+e^{-s})$, we obtain binary responses $y_{1},...,y_{n}$ by
\benn
\mathbb{P}(y_{i}=1) = p(\mbf{x}_{i}'\B),\quad 1\le i\le n,
\eenn
that is, $y_{i}$ (either 0 or 1) will be generated from a Bernoulli distribution with a probability of $p(\mbf{x}_{i}'\B)$.
\subsection{Computational aspects of the MD estimation}
We will investigate the computational aspects of the MD estimation, including the convexity of the distance function $\cL$ and the computational time required to obtain the MD estimator. Recall $\cL$ in $(\ref{eq:cL})$. Let $\B_{0}$ denote the underlying true parameter. To check the convexity of $\cL$ in the neighborhood of $\B_{0}$, we will plot its three-dimensional graph. Let $\B_{0}=(-2,1)'$ and generate $(y_{i}, \mbf{x}_{i})$, $1\le i\le 20$, accordingly. Next, create a $\B$-plane, that is, $\beta_{1}$- and $\beta_{2}$-axes, which range from -7 to 3 and from -4 to 6, respectively. Then, partition the two axes so that the length of the partitioned axis is 0.01. As a result, we obtain $10^6$ grids on the $\B$-plane. Computing $\cL$ at these grids, Figure \ref{fig:3D_loss} shows the three-dimensional graph of $\cL(\B)$ over the neighborhood of $\B_{0}$.
\begin{figure}[h]
\centering
\includegraphics[width=0.8\textwidth]{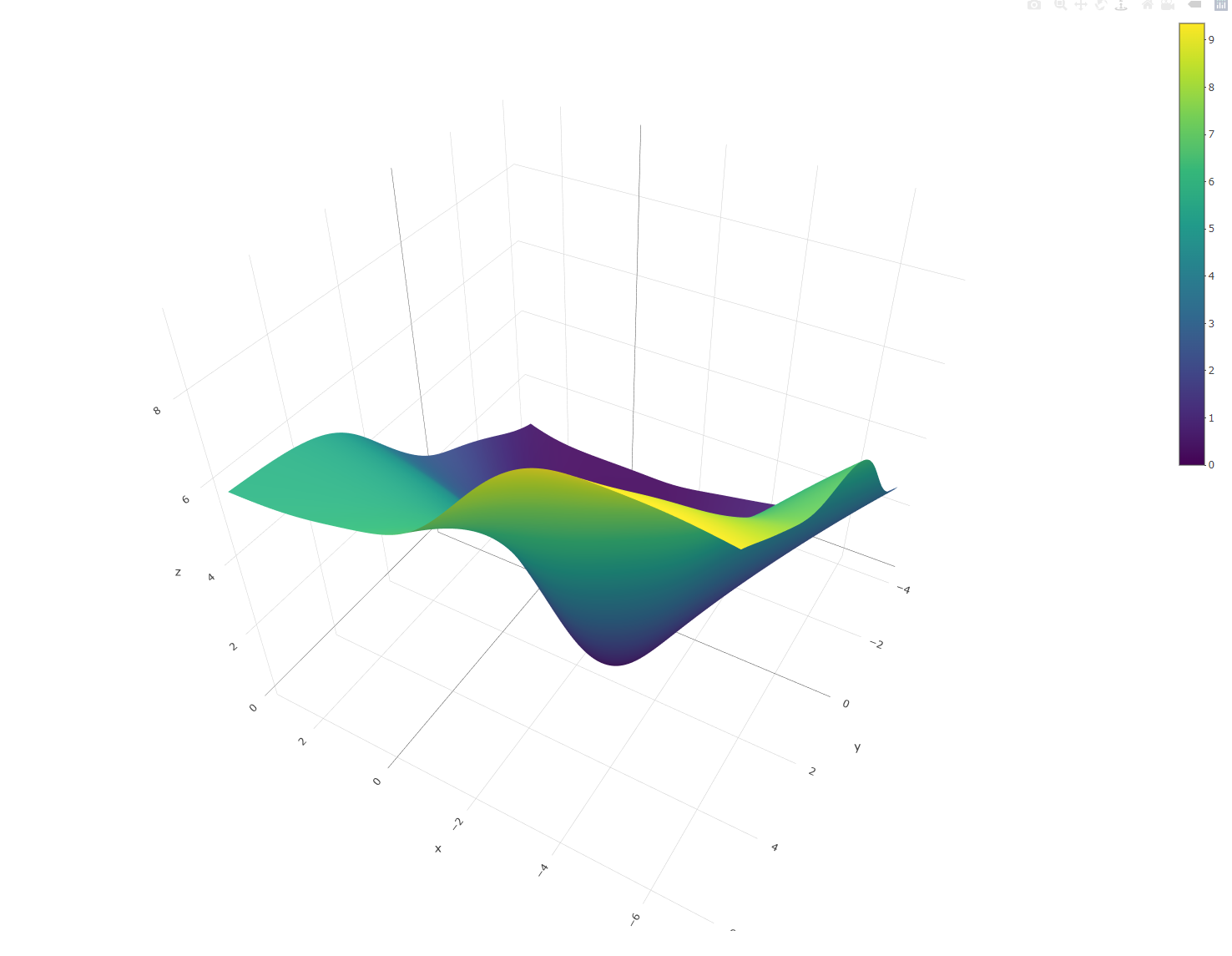}
\caption{Three-dimensional graph of $\cL(\B)$ near $\B_{0}=(-2,1)'$.}\lel{fig:3D_loss}
\end{figure}

Figure \ref{fig:3D_loss} demonstrates that $\cL$ is a convex function in the neighborhood of the true parameter, which closely accords with the ULAQ property in Section \ref{sec:ULAQ}. In addition to convexity, $\cL$ is smooth everywhere and thus differentiable. As mentioned in Remark \ref{rem:differentiability}, one caveat of the MD method in the continuous setup is its computational cost, which stems from the non-differentiability of the embedded distance function. Unlike those distance functions, the current one in our study is differentiable. Thus, when numerically solving the optimization problem in (\ref{eq:opt}) to obtain the MD estimator, we can opt for the gradient descent (GD) method. However, the GD method annexes a proviso for finding the solution: the objective function should be convex near the optimal value. Note that $\cL$ in this study is differentiable and, even convex, as shown above. Given the convexity and differentiability of $\cL$, we are motivated to use the gradient descent (GD) method to obtain the MD estimator. Recall $\mbf{S}_{n}$ in (\ref{eq:Sn}) and note that it is a gradient of the loss function, that is, $\partial \cL/\partial \B = \mbf{S}_{n}$. Thus, with starting the initial value $\B^{(0)}=(1,...,1)'\in \mathbb{R}^{J}$, we will have
\benn
\B^{(i)}=\B^{(i-1)} - lr\cdot \mbf{S}_{n},
\eenn
where the $\B^{(i-1)}$ and $\B^{(i)}$ denote the MD estimators at the $(i-1)$th and $i$th stages, respectively, and $lr$ is the learning rate used for the GD method. Then, we keep the iteration until the convergence is reached, that is,
\benn
\|\B^{(i)}-\B^{(i-1)}\|<\Delta.
\eenn
Here, $\Delta$ and $lr$ are set at 0.005 and 0.001, respectively. In the rest of this section, when computing MD estimators, we will use the GD method with the aforementioned $\Delta$ and $lr$ unless otherwise specified. \begin{figure}[h]
\centering
\begin{subfigure}[b]{0.8\textwidth}
\includegraphics[width=\textwidth]{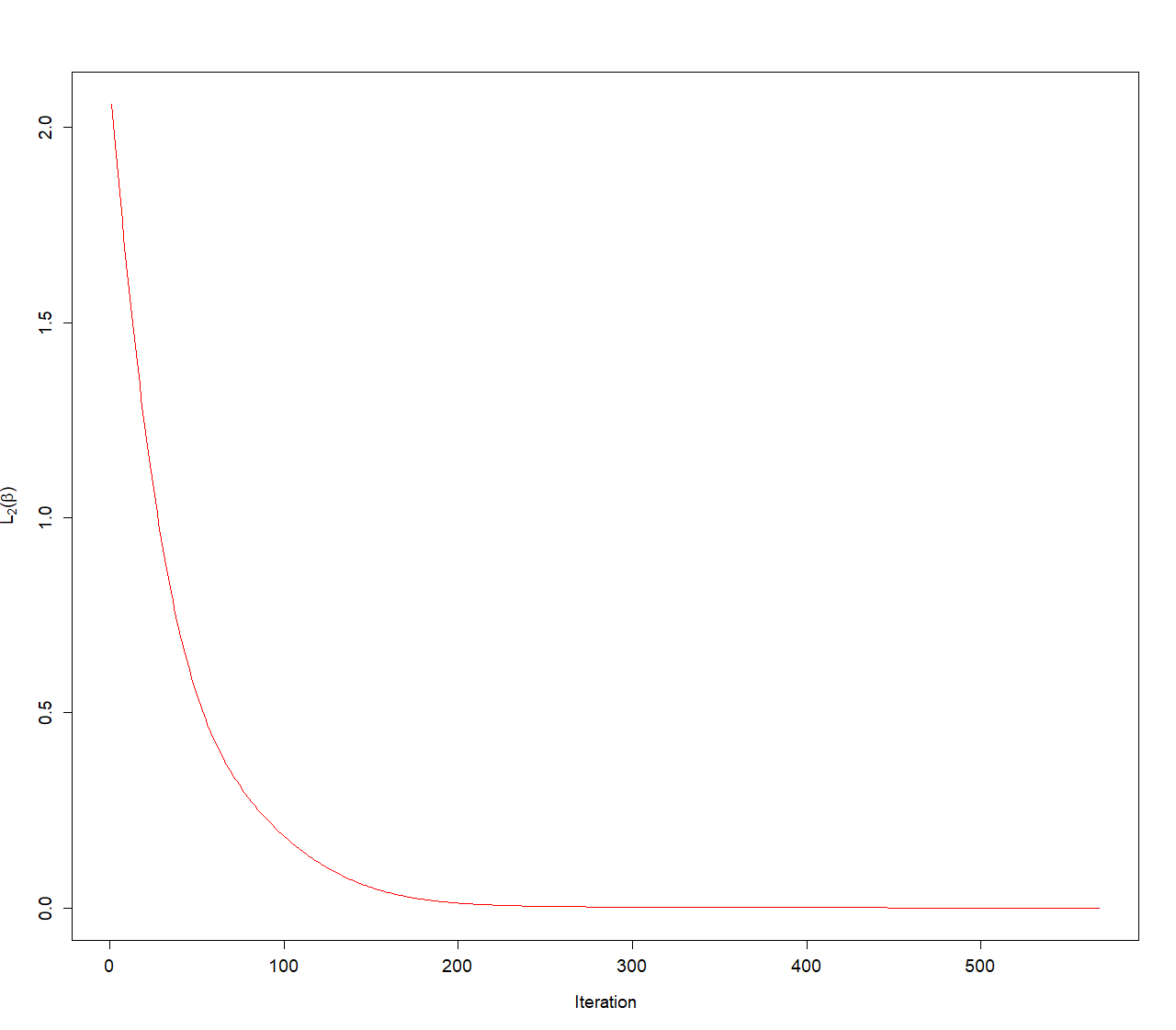}
\end{subfigure}
\caption{$\cL(\B^{(i)})$ over 1,000 iterations.}\lel{fig:graph_loss}
\end{figure}
\noi\\
As shown in Figure \ref{fig:graph_loss}, $\cL(\B^{(i)})$ decreases through the iterations, indicating that the GD method provides a solution to the optimization problem in \ref{eq:opt}, and hence, demonstrating that the convergence of the gradient descent method is actually accomplished.

Next, Table \ref{tbl:run_time} reports the computational time required to obtain the MD estimator when $n$ varies from 200 to 2,000. The second column reports the computational time itself. As shown in the table, the MD estimation yields results within 2 seconds for $n\le1,000$; even for $n = 2,000$, it takes only 3.89 seconds to complete the estimation.  Other columns in the table report the computational time divided by $n^{k}$m $k=1,2,3,4$: the computational time divided by $n$ stays around $2\times 10^{-6}$, while others decrease and never converge as $n$ increases, which empirically demonstrates that the computational time is indeed $O(n)$.
\begin{table}[h]
\centering
\begin{tabular}{|c|c|c|c|c|c|}
  \hline
  % after \\: \hline or \cline{col1-col2} \cline{col3-col4} ...
$n$ & Time & Time/$n$ & Time/$n^2$ & Time/$n^3$ & Time/$n^4$ \\
    &$(\times 10^{-3})$& $(\times 10^{-6})$ & $(\times 10^{-10})$ & $(\times 10^{-13})$ & $(\times 10^{-16})$\\
\hline
200 &0.41   &2.05   &102.5   &512.5    &2562.5\\
\hline
400 &0.8    &2      &50      &125      &312.5\\
\hline
800 &1.54   &1.925  &24.063  &30.078   &37.598\\
\hline
1000 &1.93  &1.93   &19.3    &19.3     &19.3\\
\hline
1500 &2.9   &1.933  &12.889  &8.593    &5.728\\
\hline
2000 &3.89  &1.945  &9.725   &4.862    &2.431\\
\hline
\end{tabular}
\caption{Computational time for the MD estimator with various $n$'s}.\lel{tbl:run_time}
\end{table}

\subsection{Simulation studies of bias reduction}
To assess the efficacy of bias reduction, we will compare the MD estimator without bias reduction (the original MD estimator) and the bias-reduced MD estimator described in Section \ref{sec:bias_reduction}. To compare the two estimators, this section will conduct another simulation study. Using $\B_{0}=(1.3,-2,3.5)'$ and a logistic df again, we generate a dataset as described at the beginning of this section. Using this dataset and the GD method, we obtain the original and bias-reduced MD estimators. Then, we repeat this procedure 10,000 times. After 10,000 iterations, we compute the average bias of the two MD estimators, that is, the average difference between the two MD estimators and $\B_{0}=(1.3,-2,3.5)'$. When computing the two MD estimators, we use $\bD=\bX \bA$ with $\bA=(\bX'\bX)^{-1/2}$.
\begin{table}[h]
\centering
\begin{tabular}{|c|c|c|c|c|c|c|}
  \hline
  % after \\: \hline or \cline{col1-col2} \cline{col3-col4} ...
& \multicolumn{2}{|c|}{$\b_{1}$} & \multicolumn{2}{|c|}{$\b_{2}$} & \multicolumn{2}{|c|}{$\b_{3}$} \\
\hline
$n$ & Original & Bias-reduced & Original & Bias-reduced & Original & Bias-reduced \\
\hline
20 & -0.336 & -0.173 & 1.698 & 1.660 & -1.930 & -1.596 \\
\hline
40 & -0.456 & -0.390 & 1.312 & 1.234 & -1.794 & -1.619 \\
\hline
60 & -0.490 & -0.447 & 1.153 & 1.084 & -1.698 & -1.568 \\
\hline
80 & -0.479 & -0.444 & 1.068 & 1.006 & -1.617 & -1.508 \\
\hline
\end{tabular}
\caption{Average bias of the original and bias-reduced MD estimators when $n$ varies from 20 to 80.}\lel{tbl:bias_reduction}
\end{table}
Table \ref{tbl:bias_reduction} reports the result of the bias of the original and bias-reduced MD estimators of $\beta_{i},\,i=1,2,3$. It is crystal clear that the bias-reduced MD estimators definitely yield smaller biases than the original ones, regardless of $n$ and $\beta_{i}$. In other words, compared with the original MD estimation, the MD estimation with bias reduction yields less biased estimators across all cases, indicating that bias reduction is indeed achieved.

The most significant bias reduction -- that is, 0.334 (=1.930-1.596) -- happens at $\beta_{3}$ when $n=20$, while the MD estimator of $\b_{2}$ yields the smallest reduction (0.038) when $n=20$. At first glance, the MD estimation with bias reduction appears to work well. However, the closer inspection reveals that MD estimators -- especially, $\b_{3}$ -- display rather large bias, regardless of whether or not the bias reduction is applied: the MD estimator of $\b_{1}$ yields a bias whose modulus is smaller than 0.5, while those of other MD estimators are somewhat large, ranging between 1 and 2. Considering that the bias-reduction technique can pare down the actual bias by some amount between 0.038 and 0.334, the acquisition of the original MD estimator with a not-too-large is critical. At this juncture, the merits of the MD estimation come to the fore, demonstrating its brilliance and versatility. As mentioned in Section \ref{sec:Literature}, one of the many strengths of MD estimation is that the integrating measure and weights embedded in the distance function can be tailored to the various needs of practitioners. In the setup of a continuous response variable, the MD estimation employing the Lebesgue measure and predictors for the integrating measure and weights in the distance function yields the Hodges-Lehmann type estimator, which is known to be robust; see Chapter 5.3 of Koul (2002) for more details. Motivated by the merits of the MD estimation resulting from using various weights in the distance function along with $\bD=\bX \bA$, we will try several $\bA$'s, obtain various MD estimators accordingly, and choose the optimal one. Upon trying various $\bA$'s,
empirical findings reveal that using $\bA=(\bX'\L \mathbf{P}_{n}^{-1}\L\bX)^{-1/2}$ in Remark \ref{rem:best} returns the best result, as reported in the following table.
\begin{table}[h]
\centering
\begin{tabular}{|c|c|c|c|c|c|c|}
  \hline
  % after \\: \hline or \cline{col1-col2} \cline{col3-col4} ...
& \multicolumn{2}{|c|}{$\b_{1}$} & \multicolumn{2}{|c|}{$\b_{2}$} & \multicolumn{2}{|c|}{$\b_{3}$} \\
\hline
$n$ & Original & Bias-reduced & Original & Bias-reduced & Original & Bias-reduced \\
\hline
20 & 0.018 & 0.184 & 0.636 & 0.593 & -0.757 & -0.425 \\
\hline
40 & 0.026 & 0.090 & 0.360 & 0.286 & -0.403 & -0.231 \\
\hline
60 & 0.027 & 0.069 & 0.175 & 0.105 & -0.198 & -0.067 \\
\hline
80 & 0.029 & 0.063 & 0.095 & 0.033 & -0.071 & -0.038 \\
\hline
\end{tabular}
\caption{Bias of the MD estimators before and after bias reduction.}\lel{tbl:bias_reduction2}
\end{table}

A quick glance reveals that initial MD estimators prior to the application of bias reduction display much smaller biases than those in Table \ref{tbl:bias_reduction}. The result after the application of bias reduction reported in Table \ref{tbl:bias_reduction2} is, however, slightly different from that in the previous table: the bias reduction occurs for $\b_{2}$ and $\b_{3}$ -- regardless of $n$ -- with the modulus of reduction ranging between 0.043 (corresponding to $\b_{2}$ and $n=20$) and 0.332 (corresponding to $\b_{3}$ and $n=20$), while the case of $\b_{1}$ actually shows the increase of bias for all $n$. It is presumed that the initial MD estimator of $\b_{1}$ already possesses bias small enough not to necessitate bias reduction any more: for all $n$, the modulus of the bias for $\b_{1}$ is smaller than 0.03, while those for $\b_{2}$ and $\b_{3}$ range between 0.071 and 0.757, and hence, the application of bias reduction to $\b_{1}$ gives rise to a more biased estimator. Thus, no bias reduction for $\b_{1}$ should be recommended when $\bA$ in Remark \ref{rem:best} is used for the MD estimation. In what follows, the bias-reduced MD estimator denotes the one obtained after applying bias reduction to all $\beta_{i}$'s except the coefficient of the intercept.

\subsection{MD vs.\,\,variants of GLM}
This section compares the MD estimator with the GLM estimator and its variants. As mentioned in the introduction, the MD methodology, especially with the CvM type distance, is known to yield many desirable properties. For example, Kim (2025) the MD estimation retains these properties when applied to the parameter estimation of one sample of a binomial distribution. Here, we also expect the MD estimation of the binary regression parameters to yield better (or at least similar) results than other well-celebrated estimators.

It is not difficult to see that there are many variants of a GLM in the literature. As reported in the previous section, an estimator obtained from a GLM often displayed significant bias, leading a practitioner to incorrect statistical inference. To address this issue, many researchers have proposed statistical packages for bias-reduced estimators of binary regression models. Bianco and Yohai (1998) proposed a statistical package \texttt{byglm} that provides a robust estimator when the covariates are contaminated, that is, they contain huge outliers; they demonstrated that the estimator obtained from their package remains consistent and is
asymptotically normally distributed. Based on Kosmidis and Firth (2021), Kosmidis (2025) proposed a statistical package, \texttt{brglm}, which enables practitioners to obtain a bias-reduced estimator for a binary regression model. Taking a Bayesian approach, Gelman et al.\,\,(2008) proposed a t-family of distributions -- including the Cauchy distribution -- for the prior distribution. They studied the resulting estimators for classical binary regression models. For comparison purposes, the statistical package \texttt{bayesglm}, based on Gelman et al.\,\,(2008), is used in this study. We refer to estimators obtained from \texttt{brglm}, \texttt{byglm}, and \texttt{bayesglm} as BR, BY, and Bayes estimators, respectively. Here, the MD estimator is a bias-reduced estimator introduced in the previous section. For the following analysis, we obtain the MD, BR, BY, and Bayes estimators by generating the dataset as described in the previous section; then, we compute the RMSEs and biases of those estimators for comparison.
\noi
\begin{table}[h]
\centering
\begin{tabular}{|c|c|c|c|c|c|}
  \hline
  % after \\: \hline or \cline{col1-col2} \cline{col3-col4} ...
          & $n$ & BR& BY &Bayes & MD\\
  \hline
  \hline
  \multirow{6}{*}{$\b_{1}$} &20  &0.908 (-0.327) &113.66 (36.193)  &0.846 (-0.379) &0.724 (-0.085)\\
\cline{2-6}
&40  &0.843 (-0.039) &202.536 (30.473) &0.671 (-0.186) &0.577 (-0.029)\\
\cline{2-6}
&60  &0.817 (0.026)  &94.354 (7.997)   &0.572 (-0.109) &0.513 (0.012)\\
\cline{2-6}
&80  &0.601 (-0.006) &43.569 (2.696)   &0.509 (-0.096) &0.463 (0.01)\\
\cline{2-6}
&100 &0.55 (0.017)   &38.073 (2.577)   &0.481 (-0.051) &0.438 (0.047)\\
  \hline
  \hline
  \multirow{6}{*}{$\b_{2}$} &20  &1.046 (0.482)  &172.043 (-60.253) &1.079 (0.946) &0.888 (0.61)\\
\cline{2-6}
&40  &0.879 (0.09)   &161.923 (-32.715) &0.761 (0.497) &0.525 (0.257)\\
\cline{2-6}
&60  &1.067 (-0.071) &110.106 (-11.4)   &0.673 (0.274) &0.429 (0.084)\\
\cline{2-6}
&80  &0.731 (-0.027) &56.852 (-3.811)   &0.6 (0.207)   &0.393 (0.025)\\
\cline{2-6}
&100 &0.638 (-0.025) &62.997 (-4.216)   &0.552 (0.159) &0.39 (-0.023)\\
  \hline
  \hline
\multirow{6}{*}{$\b_{3}$} &20  &1.414 (-0.834) &289.067 (105.336) &1.345 (-1.093) &0.671 (-0.415)\\
\cline{2-6}
&40  &1.327 (-0.106) &374.277 (66.159)  &1.035 (-0.521) &0.539 (-0.221)\\
\cline{2-6}
&60  &1.584 (0.044)  &165.608 (18.342)  &0.938 (-0.342) &0.479 (-0.104)\\
\cline{2-6}
&80  &1.154 (0.037)  &84.526 (6.214)    &0.881 (-0.225) &0.484 (0.007)\\
\cline{2-6}
&100 &1.003 (0.054)  &125.908 (7.994)   &0.829 (-0.146) &0.5 (0.092)\\
  \hline
\end{tabular}
\caption{Biases and SEs of the MD and other estimators. }\lel{tbl:comparison} %09072025 Lab notes.
\end{table}
Table \ref{tbl:comparison} presents the results of a comparison between the MD, BR, BY, and Bayes estimators, including the RMSEs (biases) of the estimators as $n$ varies from 20 to 100. The first six rows report the RMSEs and biases of the estimators for $\b_{1}$, while the second and third six rows report those for $\b_{2}$ and $\b_{3}$, respectively. It is worth noting several facts. First, a quick glance reveals that the BY estimator displays much larger bias, regardless of $\b_{i}$ and $n$, resulting in a larger RMSE than any other estimator. Next, all estimators display smaller bias and RMSE as $n$ increases, providing evidence of their consistency. Lastly, but most importantly, the MD estimator exhibits the best performance in terms of RMSE for all pairs of $(\b_{i}, n)$. More surprisingly, the MD estimator's superiority holds for the bias across many pairs of $(\b_{i}, n)$. For example, consider the results of biases for $\b_{1}$. For all $n$'s, the MD estimator displays the least bias for $n\le 60$ among all estimators. However, if the analysis is confined to the Bayes and MD estimators, the MD estimator shows a smaller bias than its competitor for all $n$'s. In terms of RMSE, the MD estimator outperforms the other estimators, followed by the Bayes estimator. The BY estimator performs worst due to its large bias. However, using bias as an evaluation criterion will tip the scales in favor of the BR estimator: it performs slightly better than, or at least as well as, the MD estimator.
\begin{figure}[h]
    \centering
    \begin{subfigure}[b]{0.48\textwidth}
        \includegraphics[width=\textwidth]{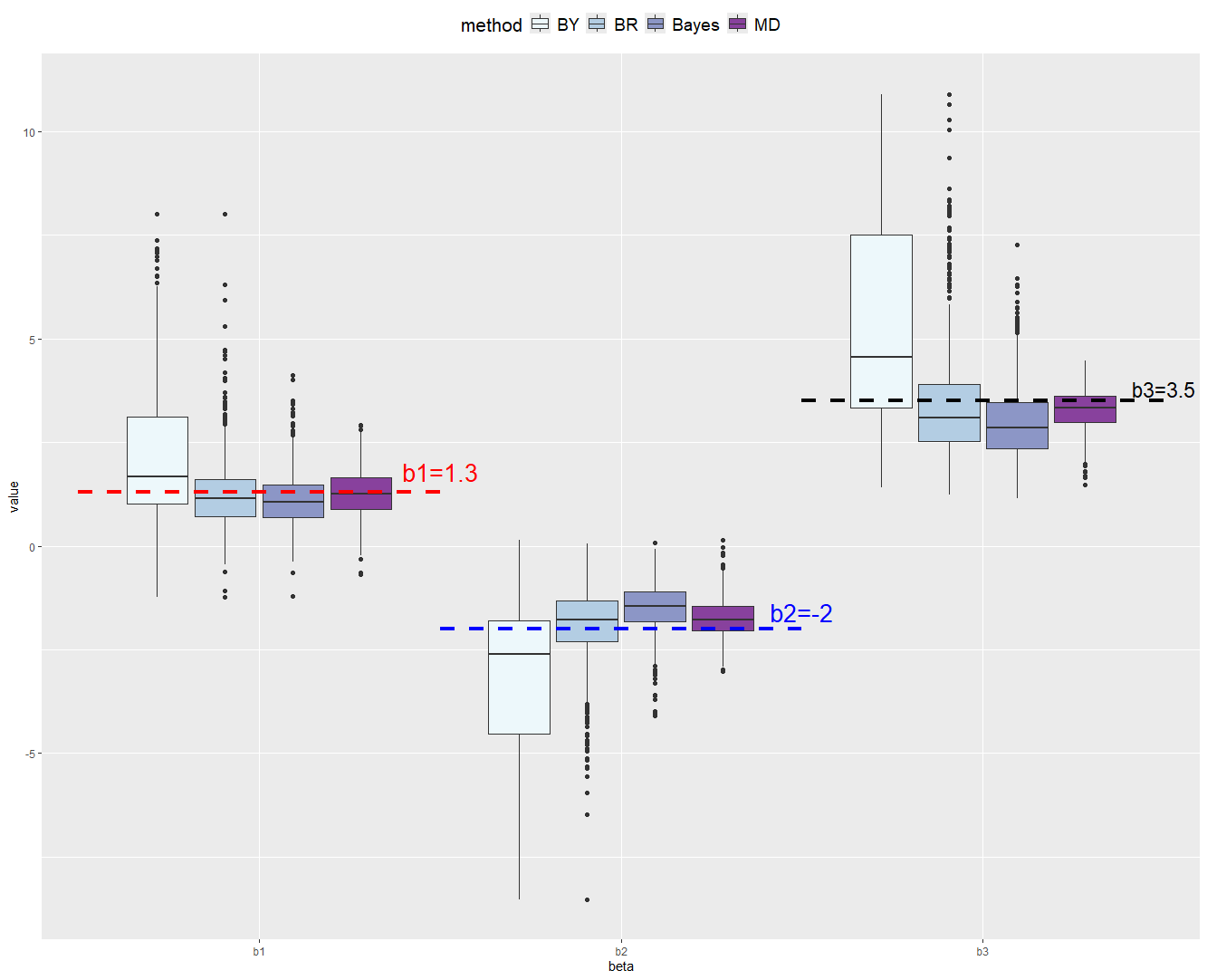}
    \end{subfigure}
    \begin{subfigure}[b]{0.48\textwidth}
        \includegraphics[width=\textwidth]{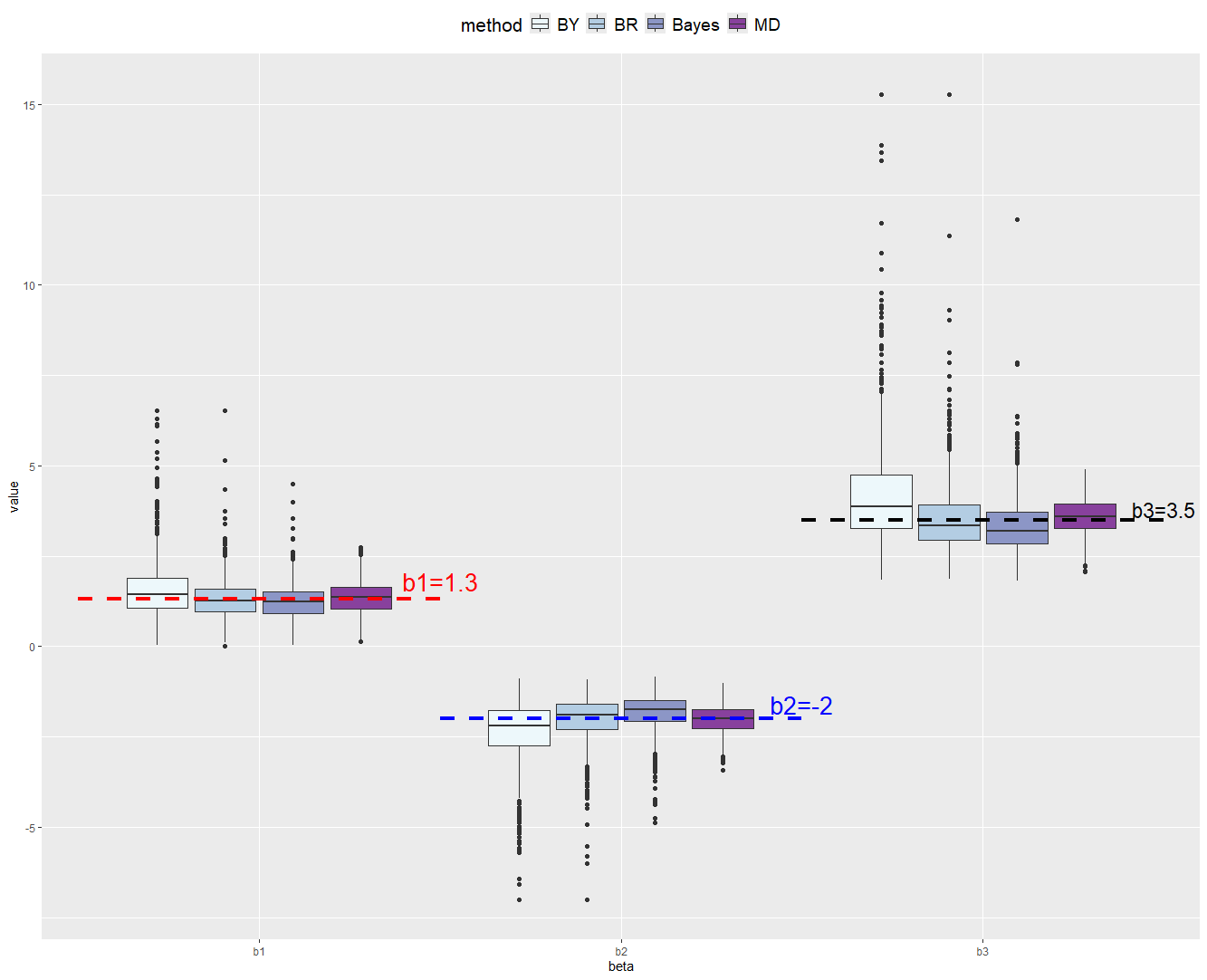}
    \end{subfigure}
    \caption{Boxplot2 of $n=40$ (left) and $n=100$ (right).}\lel{fig:bp_wo_out}
\end{figure}

The boxplots in Figure \ref{fig:bp_wo_out} show distributions of the BY, BR, Bayes, and MD estimators for $n=40$ (left panel) and $n=100$ (right panel):  the first group of boxplots shows the distributions for $\b_{1}$, while the second and third groups report those for $\b_{2}$ and $\b_{3}$, respectively. The red (1.3), blue (-2), and black (3.5) dotted lines represent the true values of $\b_{1}$, $\b_{2}$, and $\b_{3}$, respectively. First, it is worth noting that all methods show better performance -- smaller biases and SEs -- as $n$ increases from 40 to 100, which implies these estimators are consistent. Just a quick glance reveals that the BY estimators show the worst performance regardless of $n$ and $\b_{i}$, followed by the BR estimator. When only the Bayes and MD estimators are considered for comparison, with $n = 40$ (in the left panel), it is not difficult to observe that the MD estimator exhibits a smaller deviation of the median from the true value and smaller variation, and thus performs better. In the case of $n=100$ (reported in the right panel), the MD estimators still show a smaller deviation from the true values for $\b_{2}$ and $\b_{3}$, while it is quite difficult to tell for $\b_{1}$. However, it is clear to see that the MD estimators show smaller variations for all  $\b_{i}$'s, and hence, the results in Figure \ref{fig:bp_wo_out} closely accord with those in Table \ref{tbl:comparison}.

Next, we proceed to a simulation experiment to assess the predictive powers of the estimators. For the assessment, various measures will be employed: sensitivity, precision, F1, accuracy, false positive rate (FDR), and balanced accuracy (BA). In this experiment, we will generate a dataset of $n=300$. Among 300 pairs of $(y_{i}, \mbf{x}_{i})$'s, 200 will be used for training, and the rest will be allotted for validation. Using the training set, we obtain the BR, BY, Bayes, and MD estimators. Following the training session, the validation session will compute the measures above, using the estimators. Then we repeat the entire process 1,000 times and calculate averages across all measures for all estimators. Table \ref{tbl:Sim_Measure} reports the results of the experiment.
\begin{table}[h]
\centering
\begin{tabular}{|c| c| c |c |c|}
  \hline
  % after \\: \hline or \cline{col1-col2} \cline{col3-col4} ...
   Measure & BR& BY &Bayes & MD\\
  \hline
 Sensitivity &0.88  &0.884 &0.879 &0.882\\
\hline
Precision   &0.878 &0.881 &0.877 &0.882\\
\hline
F1          &0.878 &0.882 &0.877 &0.881\\
\hline
Accuracy    &0.857 &0.861 &0.856 &0.861\\
\hline
FDR         &0.122 &0.119 &0.123 &0.118\\
\hline
BA          &0.852 &0.856 &0.851 &0.856\\
 \hline
\end{tabular}
\caption{Prediction powers of the BR, BY, Bayes, and MD estimation.}\lel{tbl:Sim_Measure}
\end{table}

As reported in the table, all methods show similar performance. However, a closer look reveals that the BY and MD estimators vie for the best performance, the BR estimator follows these two estimators, and the Bayes estimator shows the worst performance. The BY estimator reports better sensitivity and F1, while the MD estimator shows better performance for precision and FDR; they tie in accuracy and BA. The BR and Bayes estimators do not show better performance than these two estimators in any single measure. Thus, it is not unnatural to conclude the simulation studies by stating that the MD method demonstrates better performance than all other methods in both estimation and prediction.
\subsection{Real examples}
\subsubsection{Vaso-constriction data set}
For the real data example, we introduce the dataset of the vaso-constriction (VC) study by Finney (1947). This dataset has been widely used by other researchers, who have incorporated it into their research works to demonstrate that their findings also successfully apply to real-world examples: see, e.g., Copas (1988) and Pregibon (1981). The VC study investigated how the volume and rate of air inspiration constrict the blood vessels in the skin; the binary response was recorded as 1 if VC occurred and 0 otherwise.
\begin{figure}[h]
    \centering
    \begin{subfigure}[b]{0.48\textwidth}
        \includegraphics[width=\textwidth]{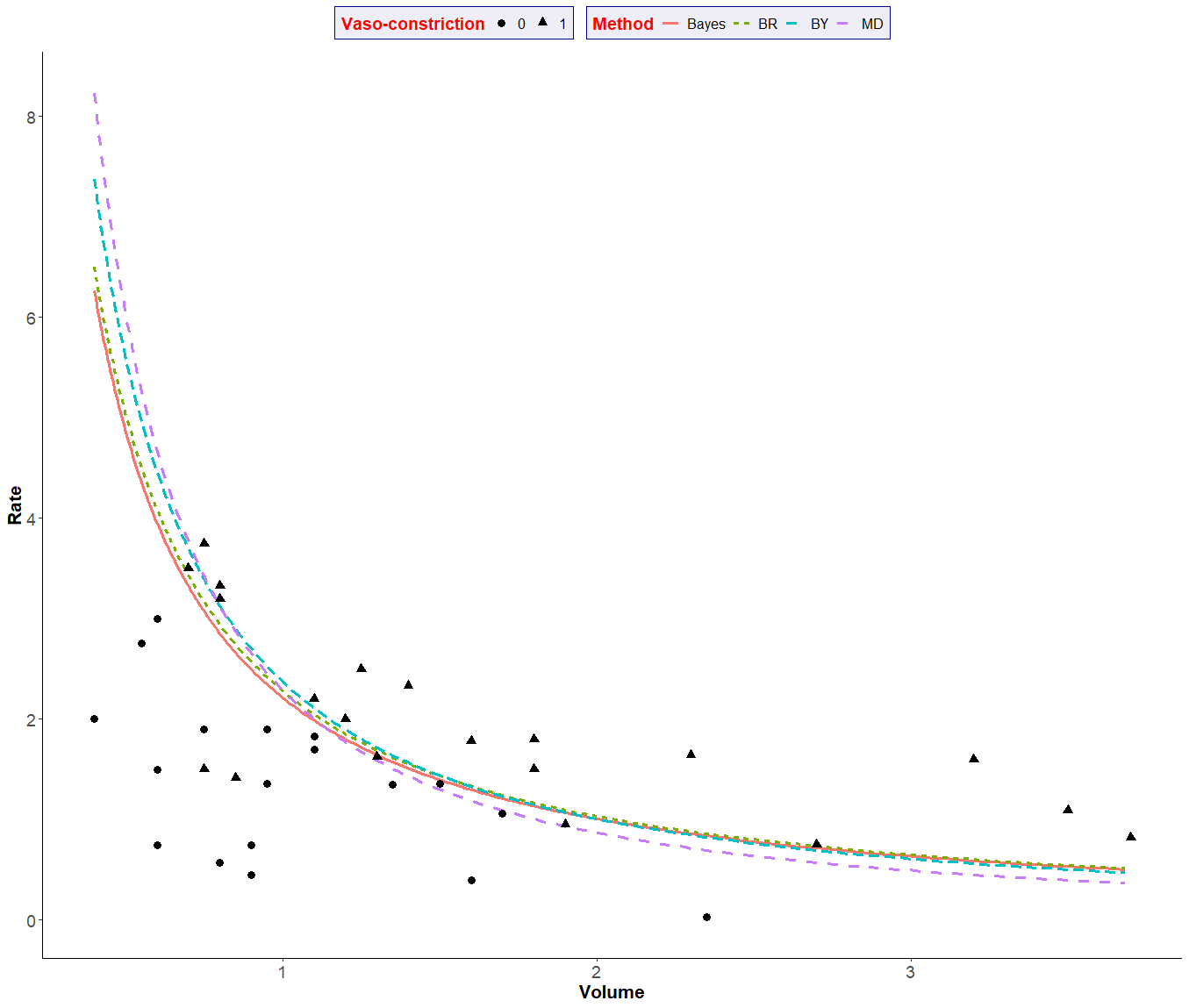}
    \end{subfigure}
    \begin{subfigure}[b]{0.48\textwidth}
        \includegraphics[width=\textwidth]{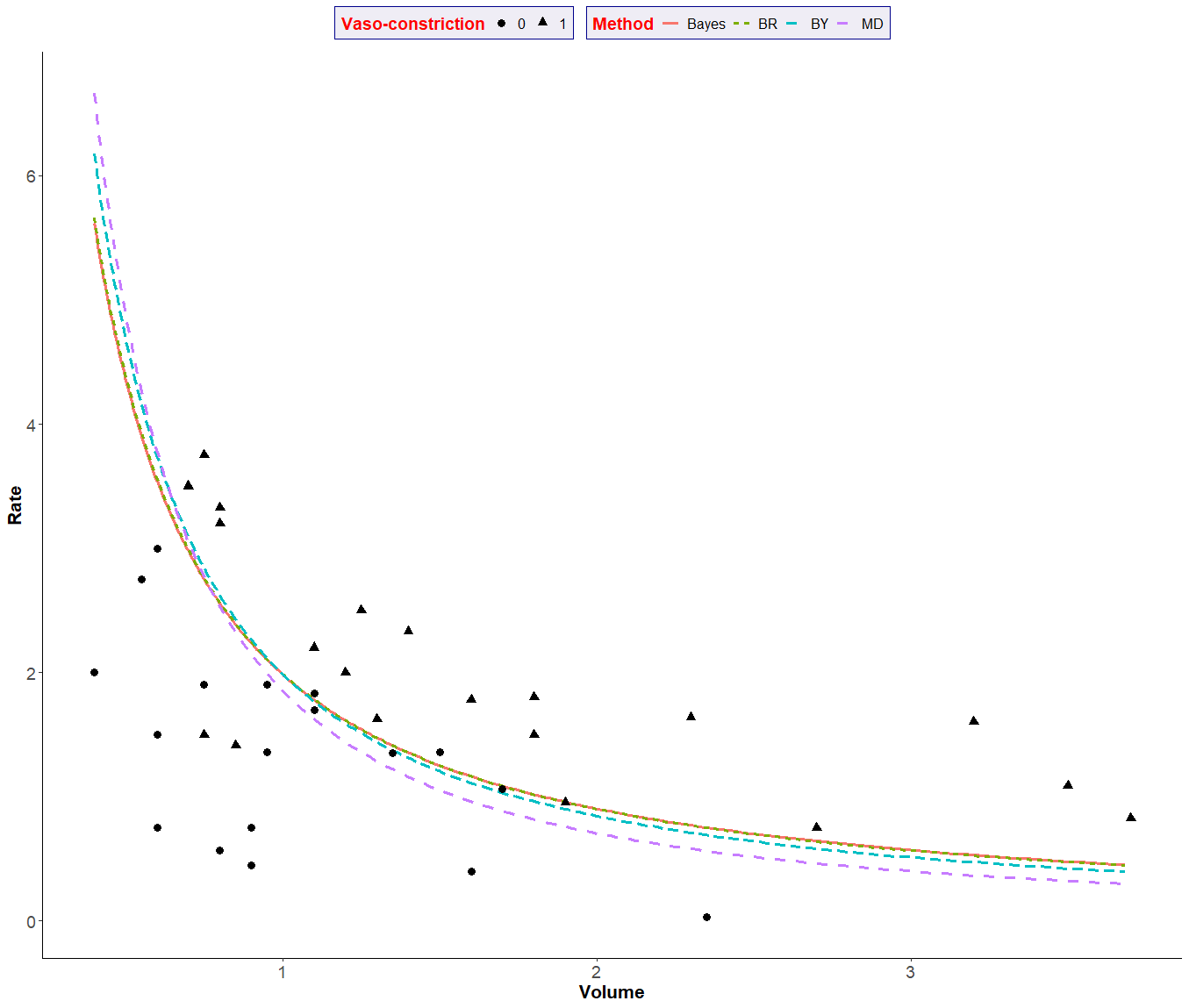}
    \end{subfigure}\caption{Scatter plot of the VC data and contour -- 0.75 (left) and 0.25 (right) -- plots of four methods with the logit link. The circle and triangle points represent the occurrence and non-occurrence, respectively, of VC. }\lel{fig:contour}
\end{figure}
Figure \ref{fig:contour} shows the scatter plot of the VC dataset and contour plots of the Bayes (orange straight), BR (light-green two-dotted), BY (skyblue long-dotted), MD (purple dotted) methods using the logistic $p$; the point marked by a triangle indicates that the VC happens, while the circular point indicates non-occurrence of VC.  The left and right subfigures present 0.75- and 0.25-contour curves, respectively, of all four methods: the 0.75-contour curve of the given method represents the fitted probability of 0.75, implying that the method predicts any point above the curve to have a probability of the VC status greater than 0.75. As described in the figure, all methods exhibit similar 0.75- and 0.25-contour curves.
\subsubsection{Urinary incontinence data set}
For the second real example, we utilize the dataset from the urinary incontinence (UI) study by Potter (2005); for more details, see Table 1 therein. When the UI study applied treatments to patients, the binary response was recorded as 1 if the patients could control their bowel and bladder, and 0 otherwise. Potter (2005) investigated the relationship between the binary response and three urinary tract variables ($X_{1}, X_{2}$, and $X_{3}$) using a binary regression.
\begin{figure}[h]
\centering
    \begin{subfigure}[b]{0.48\textwidth}
        \includegraphics[width=\textwidth]{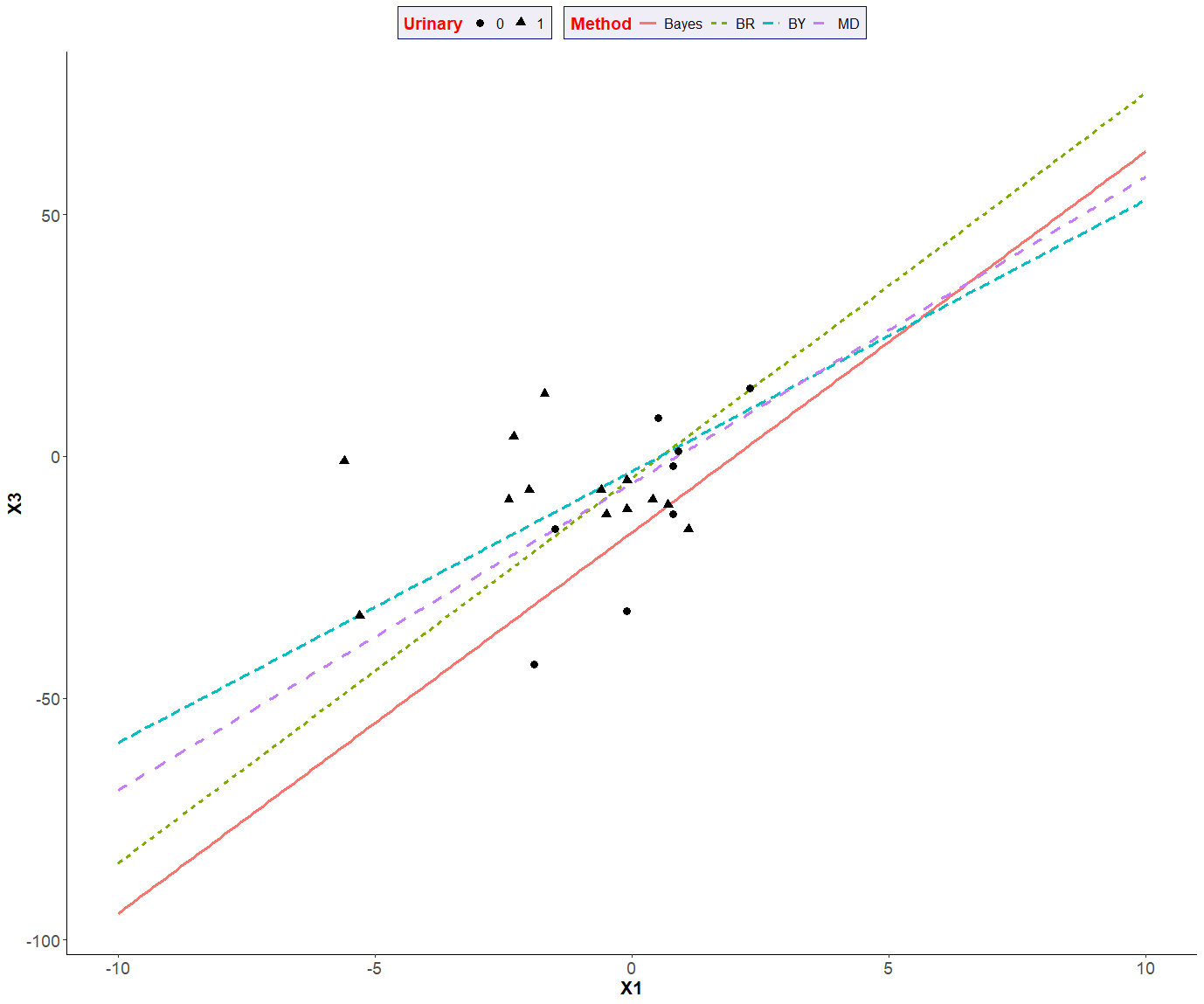}
    \end{subfigure}
    \begin{subfigure}[b]{0.48\textwidth}
        \includegraphics[width=\textwidth]{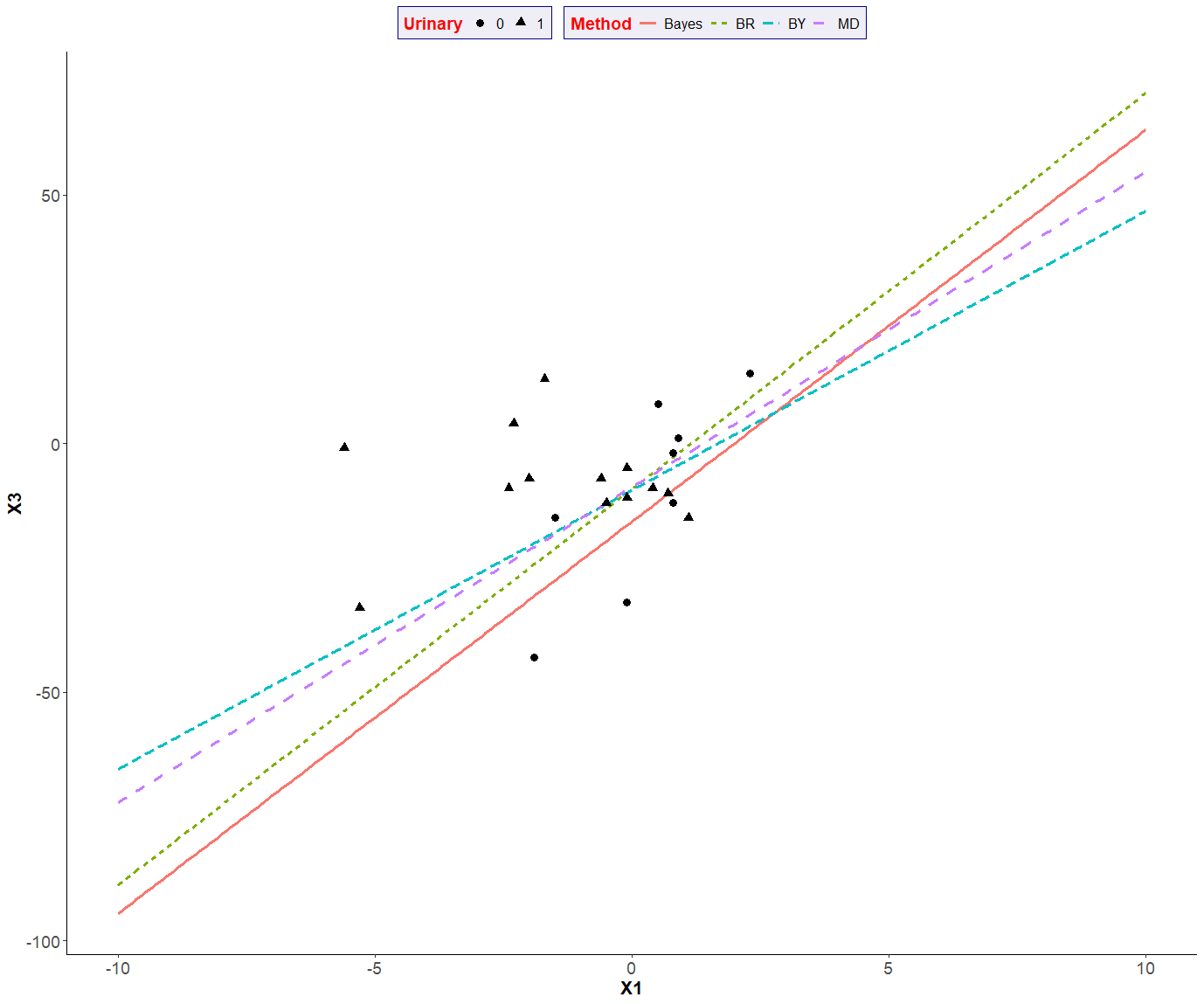}
    \end{subfigure}
\caption{Scatter plot of the UI data set and contour -- 0.75 (left) and 0.25 (right) -- plots of four methods with the logit link. The circle and triangle points represent the urinary incontinence and continence, respectively.}\lel{fig:UI_contour}
\end{figure}
Figure \ref{fig:UI_contour} shows 0.75- and 0.25-contour plots obtained from application of the four methods to the UI dataset. To generate two-dimensional contour plots, we fix $X_{2}$ at -3 and use the logit link function again. Unlike the VC dataset, glaringly obvious differences exist between the methods. For example, both 0.75- and 0.25-contour lines of the Bayes method are parallel to those of the BY method, while displaying larger slopes than those of the BY and MD methods.

Next, we apply another binary regression with the Cauchit link function to the UI dataset and see whether using it yields any different results compared with those obtained from using the logit link function. When comparing contour plots obtained with the logit and Cauchit link functions, we report the comparison for the MD method only, as it exhibits the most pronounced difference between the two link functions.
\begin{figure}[h]
\centering
\includegraphics[width=0.6\textwidth]{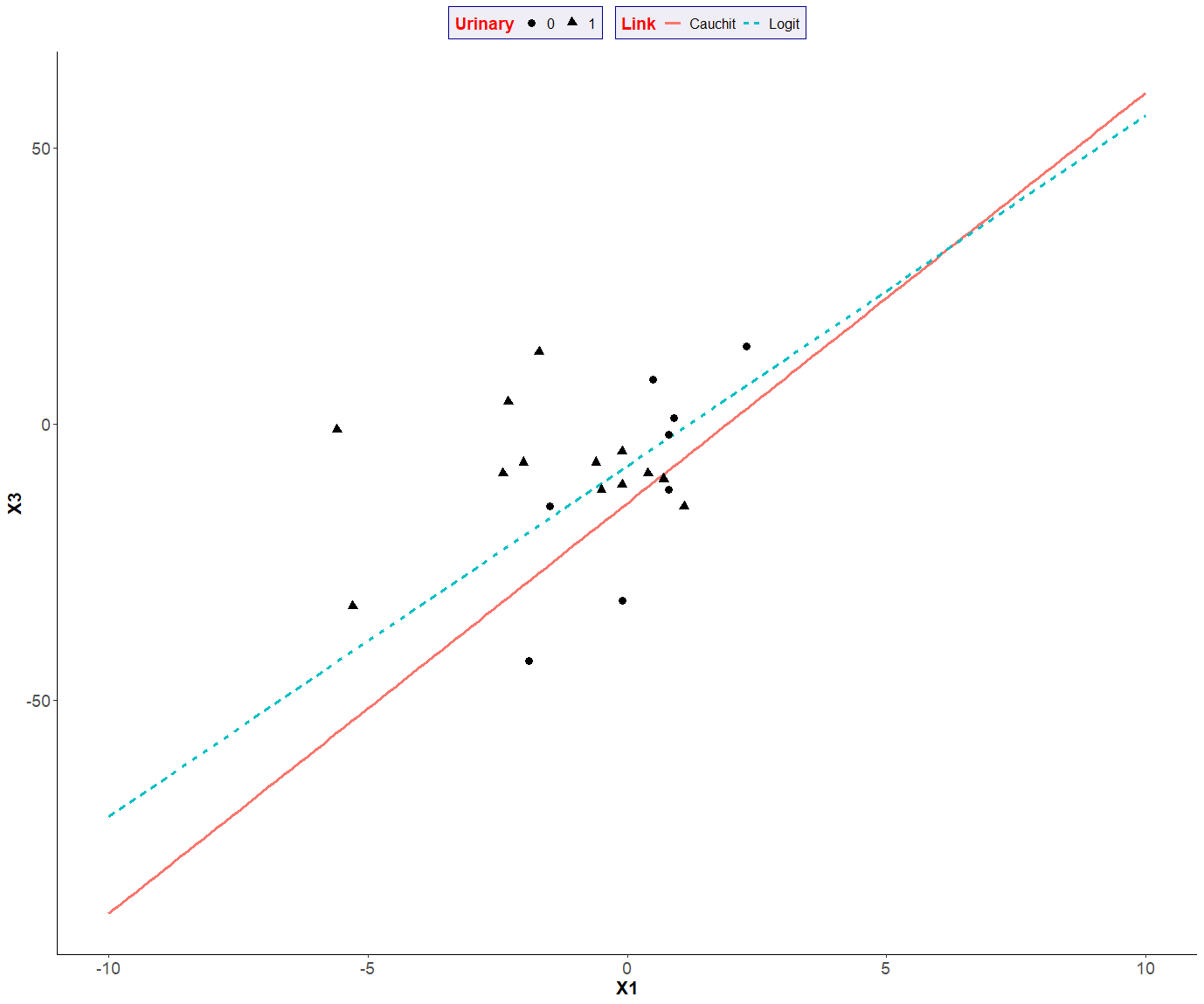}
\caption{Scatter plot of the UI dataset and 0.75-contour plots of the MD method with the logit (blue) and Cauchit (red) link functions.}\label{fig:UI_contour_cauchit}
\end{figure}
Figure \ref{fig:UI_contour_cauchit} shows 0.75-contour plots corresponding to the logit (blue) and Cauchit (red) link functions. As illustrated in the figure, employing a different link function results in a significant discrepancy, even with the same MD method. Furthermore, the 0.75-contour line obtained from the Cauchit link function contains more urinary continence cases, which implies the MD method with the Cauchit link function provides a better model fitting than the logit counterpart.

The visual comparison illustrated in the previous figures has a dearth of information to claim superiority of a particular method or link function. For a more concrete comparison, we therefore obtain numerical measures that show the diagnostic accuracy of the chosen method and link function. To this end, five measures -- sensitivity, precision, F1 score, accuracy, and false discovery rate (FDR) -- will be used for comparison purposes. When computing these measures, we will employ cross-validation. Among various cross-validation methods, we will use leave-one-out cross-validation (LOOCV), as the UI dataset comprises only 21 observations. More precisely, after removing the first observation for validation purposes and using the remaining dataset, we obtain the BR, BY, Bayes, and MD estimators. Using the estimators, obtain the predicted probability of the occurrence of UI and save it. We then repeat the procedure 21 times, reserving the i-th observation for the i-th validation each time. Consequently, having 21 pairs of ground truth labels and predictions, we proceed to find the five aforementioned measures. Table \ref{tbl:urinary} reports the result of the LOOCV.

The table below reports the results of the LOOCV analysis when the logit and Cauchit link functions are employed. Numerical figures in the table represent measures of the corresponding methods using the logit link function, while those in parentheses represent measures of the same methods using the Cauchit link function. For example, the BR method, using the logit and Cauchit link functions, yields 0.75 and 0.778 for the accuracy, respectively, while the Bayes method, using the same link functions, yields 0.74 and 0.735 for the same measure, respectively. Since the BY method doesn't provide an option for the Cauchit link function, ``Not Available" (NA) is reported for the Cauchit case.
\begin{table}[h]
\centering
\begin{tabular}{|c| c| c |c |c|}
  \hline
  % after \\: \hline or \cline{col1-col2} \cline{col3-col4} ...
   Measure & BR& BY &Bayes & MD\\
  \hline
Sensitivity &0.818 (0.826) &0.846 (NA) &0.796 (0.786) &0.824 (0.844)\\
\hline
Precision   &0.792 (0.819) &0.786 (NA) &0.792 (0.793) &0.877 (0.906)\\
\hline
F1          &0.803 (0.822) &0.815 (NA) &0.791 (0.786) &0.848 (0.873)\\
\hline
Accuracy    &0.75 (0.778)  &0.762 (NA) &0.74 (0.735)  &0.816 (0.849)\\
\hline
FDR         &0.208 (0.181) &0.214 (NA) &0.208 (0.207) &0.123 (0.094)\\
\hline
\end{tabular}
\caption{Prediction powers of the BR, BY, Bayes, and MD estimations with the logit and Cauchit link functions for the urinary data set.}\lel{tbl:urinary}
\end{table}
The results reported in Table \ref{tbl:urinary} demonstrate that the MD method outperforms all other methods, regardless of which measures and link functions are employed. When the logit link function is used, the MD method exhibits the best performance, followed by the BR and BY methods; overall, the Bayes method performs worst. Using the Cauchit link function does not change the balance tipped in favor of the MD method: the MD method still outperforms other methods, followed by the BR, while the Bayes method again shows the worst performance. When comparison is made between the logit and Cauchit link functions for a given method, an interesting fact emerges. More precisely, both MD and BR methods exhibit more optimal results for all measures (i.e., smaller value for FDR and larger value for all other measures) when the Cauchit link function is used, while the Bayes method exhibits the opposite outcome. %Additional investigation reveals more remarkable facts.

In conclusion, the MD method outperforms other methods for both logit and Cauchit link functions. However, the extent of superiority of the MD method to the other methods becomes more prominent when the Cauchit link function is used; e.g., the MD and Bayes methods corresponding to the logit link function report 0.824 and 0.796 for the sensitivity, respectively, while the Cauchit counterparts report 0.844 and 0.786, respectively, thereby showing a larger difference (0.028 vs. 0.058). Even though the extent of the superiority of the Cauchit link function to the logit link function is diminished a bit, the same fact holds for the MD vs. BR methods: for the sensitivity, two methods exhibit the difference of 0.006 (=0.824-0.818) and 0.018 (=0.844-0.826) for the logit and Cauchit link functions, respectively. Note that the results reported in Table \ref{tbl:urinary} closely accord with those described in Figure \ref{fig:UI_contour_cauchit}, where the 0.75-contour curve corresponding to the Cauchit link function exhibits a better model fitting than that corresponding to the logit link function.

\section{Conclusion}\lel{Sec:conclusion}
This study extends the application of the $L_{2}$ optimization using the CvM type distance function from one sample setup of Kim (2025) to a binary regression setup, proposes the MD estimator of the regression parameters, and investigates its asymptotic properties. Furthermore, by emulating methodologies suggested in the literature on the GLM estimation, this study proposes bias reduction of the MD estimator and demonstrates that the performance of the MD method can be improved to some degree. Simulation studies and real examples demonstrate that the proposed method compares favorably with other well-celebrated methods. The current research can serve as a benchmark for further applications of MD estimation in other regression analyses (multinomial logistic regression, Poisson regression, etc.), which will form the basis of future research.

\noi\\

\edt

\begin{table}[h]
\centering
\begin{tabular}{|c| c| c |c |c|}
  \hline
  % after \\: \hline or \cline{col1-col2} \cline{col3-col4} ...
   Measure & BR& BY &Bayes & MD\\
  \hline
Sensitivity &0.8   &0.65  &0.85  &0.85\\
\hline
Precision   &0.8   &0.812 &0.708 &0.739\\
\hline
F1          &0.8   &0.722 &0.773 &0.791\\
\hline
Accuracy    &0.795 &0.744 &0.744 &0.769\\
\hline
FDR         &0.2   &0.188 &0.292 &0.261\\
\hline
BA          &0.795 &0.746 &0.741 &0.767\\
  \hline
\end{tabular}
\caption{Biases and SEs of the MD and other estimators.}\lel{tbl:CV_Logit}
\end{table}

\begin{table}[h]
\centering
\begin{tabular}{|c| c| c |c |c|}
  \hline
  % after \\: \hline or \cline{col1-col2} \cline{col3-col4} ...
   Measure & BR& BY &Bayes & MD\\
  \hline
Sensitivity &0.7   &0.65  &0.4   &0.75\\
\hline
Precision   &0.7   &0.812 &0.571 &0.75\\
\hline
F1          &0.7   &0.722 &0.471 &0.75\\
\hline
Accuracy    &0.692 &0.744 &0.538 &0.744\\
\hline
FDR         &0.3   &0.188 &0.429 &0.25\\
\hline
BA          &0.692 &0.746 &0.542 &0.743\\
\hline
\end{tabular}
\caption{Biases and SEs of the MD and other estimators.}\lel{tbl:CV_Probit}
\end{table}

\begin{table}[h]
\centering
\begin{tabular}{|c| c| c |c |c|}
  \hline
  % after \\: \hline or \cline{col1-col2} \cline{col3-col4} ...
   Measure & BR& BY &Bayes & MD\\
  \hline
Sensitivity &0.6   &0.8   &0.65  &0.55\\
\hline
Precision   &0.667 &0.842 &0.929 &0.786\\
\hline
F1          &0.632 &0.821 &0.765 &0.647\\
\hline
Accuracy    &0.641 &0.821 &0.795 &0.692\\
\hline
FDR         &0.333 &0.158 &0.071 &0.214\\
\hline
BA          &0.642 &0.821 &0.799 &0.696\\
  \hline
\end{tabular}
\caption{Biases and SEs of the MD and other estimators.}\lel{tbl:CV_LogLog}
\end{table}

\begin{table}[h]
\centering
\begin{tabular}{|c| c| c |c |c|}
  \hline
  % after \\: \hline or \cline{col1-col2} \cline{col3-col4} ...
   Measure & BR& BY &Bayes & MD\\
  \hline
Sensitivity &0.7   &0.7   &0.75  &0.8\\
\hline
Precision   &0.875 &0.875 &0.882 &0.64\\
\hline
F1          &0.778 &0.778 &0.811 &0.711\\
\hline
Accuracy    &0.795 &0.795 &0.821 &0.667\\
\hline
FDR         &0.125 &0.125 &0.118 &0.36\\
\hline
BA          &0.797 &0.797 &0.822 &0.663\\
\hline
\end{tabular}
\caption{Biases and SEs of the MD and other estimators.}\lel{tbl:CV_Cauchit}
\end{table}

felicitous, fetching, weave, abate, across the board, act in concert with, several factors are adduced to explain..., adventitious, affirmative, affirm, agent

outlandish, outrageous, paltry